\documentclass{article}
\usepackage{amstex}

\title{Automorphic hyperfunctions and period functions}
\author{Roelof W.~Bruggeman}

\newenvironment{stel}[1]{\refstepcounter{subsection}%
\trivlist\item[]{{\bf\thesubsection}\ \,{\bf #1. }}%
\begin{sl}\ignorespaces}{\end{sl}\endtrivlist}

\def\ntstamount{2pt plus1pt minus1pt} 
\newcommand{\ntst}[1]{\vspace\ntstamount\par\noindent%
\refstepcounter{subsection}{\bf\thesubsection}%
{\def\testempt{#1}
\ifx\testempt\empt . \
\else{%
\ \,{\sl #1}\/
}\fi}%
\ignorespaces}

\newcommand{\bw}[1]{\par\hspace{-\parindent}{\sl
#1. }\ignorespaces}

\let\vgl=\eqref
\newcommand\nwsect[1]{\setcounter{equation}{0}\section{#1}}

\usepackage{amssymb}
\renewcommand\setminus{\smallsetminus}

\newcommand\hv{\mathfrak{H}}
\newcommand\glie{{\mathfrak g}}

\newcommand\CC{\mathbb C}
\newcommand\NN{\mathbb N}
\newcommand\PP{\mathbb P}
\newcommand\QQ{\mathbb Q}
\newcommand\RR{\mathbb R}
\newcommand\ZZ{\mathbb Z}

\newcommand\EE{{\mathbf E}}
\newcommand\HH{{\mathbf H}}
\newcommand\WW{{\mathbf W}}
\newcommand\XX{{\mathbf X}}
\newcommand\YY{{\mathbf Y}}

\newcommand\al{\alpha}
\newcommand\bt{\beta}
\newcommand\G{\Gamma}
\newcommand\g{\gamma}
\newcommand\dt{\delta}
\newcommand\e{\varepsilon}
\newcommand\z{\zeta}
\renewcommand\th{\theta}
\renewcommand\k{\kappa}
\newcommand\Ld{\Lambda}
\newcommand\ld{\lambda}
\newcommand\m{\mu}
\newcommand\n{\nu}
\newcommand\X{\Xi}
\newcommand\x{\xi}
\newcommand\p{\pi}
\newcommand\R{{\sf P}}
\renewcommand\r{\rho}
\newcommand\s{\sigma}
\newcommand\ups{\upsilon}
\newcommand\ph{\varphi}

\newcommand\Ps{\Psi}
\newcommand\ps{\psi}
\newcommand\om{\omega}

\newcommand\Akr{{\mathcal A}}
\newcommand\Bkr{{\mathcal B}}
\newcommand\Fkr{{\mathcal F}}
\newcommand\Hkr{{\mathcal H}}
\newcommand\Lkr{{\mathcal L}}
\newcommand\Okr{{\mathcal O}}

\newcommand\re{\operatorname{Re}}
\newcommand\im{\operatorname{Im}}

\newcommand\sign{\operatorname{sign}}
\newcommand\res{\operatorname*{res}}
\newcommand\supp{\operatorname{Supp}}

\newcommand\Gf{\operatorname{\Gamma}}
\newcommand\oh{O}

\newcommand\isdef{\mathrel{:\mskip2mu=}}
\newcommand\divides{\mathrel{|}}

\newcommand\vz[1]{\mathchoice{\left\{ #1
\right\}}{\left\{ #1 \right\}}{\{ #1 \}}{\{ #1 \}}}
\newcommand\vzm[2]{\mathchoice{\left\{\, #1 : #2
\,\right\}}{\{\, #1 :\allowbreak #2 \,\}}{\{ #1
:\allowbreak #2 \}}{\{ #1 :\allowbreak #2 \}}}
\newcommand\matc[4]{\left({#1\atop#3}{#2\atop#4}\right)}
\newcommand\matr[4]{\left({\hfill#1\atop\hfill#3}
{\hfill#2\atop\hfill#4}\right)}
\newcommand\widearray[1]{\renewcommand\arraystretch{1.4}
\begin{array}{#1}}

\newcommand\SL{{\mathrm{SL}}}
\newcommand\PSL{{\mathrm{PSL}}}
\newcommand\PGL{{\mathrm{PGL}}}

\newcommand\bhv{\hv^+}
\newcommand\ohv{\hv^-}
\newcommand\Gmod{\G_{\rm mod}}
\newcommand\Psimod{\Ps_{\rm mod}}
\newcommand\Fkrmod{\Fkr_{\rm mod}}
\newcommand\prc{\PP^1_\CC}

\newcommand\prq{\PP^1_\QQ}
\newcommand\spgl{\operatorname{\mathfrak j}}

\newcommand\pr{\operatorname{pr}}
\newcommand\pvint{\sideset{\rm pv}{\negthickspace}\int}

\renewcommand\mapsto{\mathrel{\mapstochar \mskip1mu
\rightarrow}}

\date{}
\begin{document}
\maketitle

\nwsect{Introduction}\label{secintro}
Lewis, \cite{Lw}, has given a relation between Maass
forms and period functions. This paper investigates
that relation by means of the invariant hyperfunctions
attached to automorphic forms.

\ntst{Maass forms.}\label{mafo}A cuspidal Maass form is
a function on the upper half plane $\bhv=\vzm{z\in\CC}
{\im z>0}$ that satisfies $u(-1/z)=u(z)$ and has an
expansion
\begin{equation}
\label{Maassform}
u(z) = \sum_{n\neq0} a_n W_{0,s-1/2}(4\p|n|y)e^{2\p
inx}.
\end{equation}
We write $x=\re z$ and $y=\im z$ for $z\in \bhv$, and
use the Whittaker function
$W_{\cdot,\cdot}(\,\cdot\,)$, see, e.g.,
\cite{Sl},~1.7. One can express $W_{0,\cdot}$ in terms
of a modified Bessel function: $W_{0,\m}(y) =
\sqrt{y/\p} K_\m(y/2)$.

These Maass forms occur as eigenfunctions in the
spectral decomposition of the Laplacian in $L^2\left(
\Gmod\backslash \bhv,\frac{dx\,dy}{y^2}\right)$, with
$\Gmod\isdef\PSL_2(\ZZ)$. The eigenvalue is
$s\left(1-s\right)$. For any given~$s$ the space of
such Maass forms has finite dimension. The dimension
is non-zero only for a infinite discrete set of points
on the line $\re s=\frac12$. For more information
concerning Maass forms see, e.g., \cite{Te}, \S3.5--6.

One calls a Maass form {\sl even}, respectively {\sl
odd}, if $u(-\bar z)= u(z)$, respectively $u(-\bar
z)=-u(z)$. In terms of the Fourier coefficients this
amounts to $a_{-n} = a_n$, respectively $a_{-n}=-a_n$.

Although spectral theory states that cuspidal Maass
forms exist, none of them is explicitly known. There
are computational results, see, e.g., \cite{St},
\cite{GS84}, and \cite{He92}.

\ntst{Period functions.}In \cite{Lw} and~\cite{LZ},
Lewis and Zagier show that there is a bijective linear
map from the space of cuspidal Maass forms of
weight~$0$ for a fixed value of~$s$ to the space of
holomorphic functions $\ps:\CC\setminus(-\infty,0]
\rightarrow\CC$ that satisfy
\begin{equation}
\label{MFE}
\ps(z)-\ps(z+1) = \left(z+1\right)^{-2s}
\ps\left(\frac{z}{z+1} \right),
\end{equation}
and $\ps(1)=0$,
$\lim_{z\rightarrow\infty,z\in\RR}\ps(z)=0$.
In~\cite{Lw} this bijection is given in terms of a
sequence of integral transforms; the approach
in~\cite{LZ} uses the $L$-series attached to the Maass
form. In~\cite{Za} Zagier gives indications that the
function~$\ps$ generalizes the period polynomial
associated to holomorphic cusp forms. So the name {\sl
period function}\/ is appropriate.

For even, respectively odd Maass forms, Lewis and
Zagier show that $\ps(z) = \pm
z^{-2s}\ps\left(\frac1z\right)$. Under this
assumption, equation \vgl{MFE} is equivalent to
\begin{equation}
\label{SFE}
\ps(z)-\ps(z+1) = \pm z^{-2s}\ps\left(1+\frac1z\right).
\end{equation}
Applying \vgl{SFE} to $z$, $z+1$, $z+2$, $\ldots$, and
using the behavior as $z\rightarrow\infty$, Lewis and
Zagier find
\begin{equation}
\label{efTO}
\pm \ps(z) = \sum_{n=0}^\infty (z+n)^{-2s}
\ps\left(1+\frac1{z+n} \right).
\end{equation}
This means that $\ps$ corresponds to an eigenfunction
of the {\sl transfer operator} of Mayer,~\cite{May}.
Theorem~2 in~\cite{May} shows that the
$\pm1$-eigenvectors of the transfer operator are
closely related to the zeros of the Selberg zeta
function. Lewis remarks that not only the cuspidal
Maass forms, but also some Eisnestein series should
yield eigenvectors. Zagier, \cite{Za}, indicates how
to obtain period functions by meromorphic continuation
of a partial Eisenstein series. For
$\z\left(2s\right)=0$ these functions are
eigenfunctions of the transfer operator.

\ntst{Boundary form.}Lewis, \cite{Lw}, \S6 (c), gives
formal computations with the {\sl boundary form}\/
associated to a Maass form as a motivation for his
method. The present paper arose from the wish to make
this precise, and to understand the map $u\mapsto\ps$
from Maass forms to period functions in terms of the
boundary form.

The boundary forms that we use are hyperfunctions on
the boundary of the upper half plane. These
hyperfunctions are $\SL_2(\ZZ)$-invariant vectors in a
principal series representation. Actually, the
hyperfunctions related to the most interesting
automorphic forms are distributions; that aspect we do
not discuss in this paper.

The hyperfunction point of view turns out to give two
interpretations of the period function~$\ps$. The
first one, in Theorem~\ref{authftopsi}, arises
naturally when describing any hyperfunction associated
to a modular form (even for exponentially increasing
modular forms). The second interpretation is more
complicated. We shall consider a type of parabolic
cohomology with values in the hyperfunctions. In
Proposition~\ref{coh} we show that the hyperfunctions
associated to a class of modular forms (containing the
cuspidal Maass forms and some Eisenstein series)
correspond to classes in these cohomology groups. In
Section~\ref{secrho} we use a map from hyperfunctions
to holomorphic functions on the upper half plane to
arrive at cohomology classes with holomorphic
functions as values. Such a class is determined by one
function that turns out to satisfy~\vgl{MFE}, but with
$s$ replaced by $1-s$.

It should be emphasized that we do not recover all
results of Lewis and Zagier. We do not prove that each
period function satisfying~\vgl{MFE} with the
prescribed behavior at~$1$ and~$\infty$ comes from a
Maass form.

\ntst{Geodesic decomposition.}To obtain the second
interpretation of the period function, we use what we
call {\sl $\G$-decompositions}\/ of hyperfunctions.

Let the boundary of the upper half plane be written as
a finite union $\bigcup_{j=1}^n I_j$, with closed
intervals $I_j$ that intersect each other only in
their end points, and where the end points are cusps.
Any hyperfunctions~$\al$ on the boundary of~$\bhv$ can
be written as a sum $\sum_{j=1}^n \al_j$, such that
the support of the hyperfunction~$\al_j$ is contained
in~$I_j$. There are many possibilities to arrange
this. For the hyperfunctions associated to cuspidal
Maass forms, holomorphic modular cusp forms, and some
Eisenstein series, this can be done in a neat way,
which we shall call the geodesic decomposition. We
shall show in~\ref{perfholcf} that the cocycles
attached to holomorphic cusp forms can be derived from
this decomposition.

\ntst{Overview.}From the representational point of view
it is more convenient to work with modular forms on
the group $\PSL_2(\RR)$ than on the upper half
plane~$\bhv$. This step is carried out in
Section~\ref{secautforms}. Actually, we do not
restrict ourselves to the modular group, but work with
a general cofinite discrete subgroup~$\G$, that is
required to possess cusps. Section~\ref{sechf}
discusses some properties of hyperfunctions.
Section~\ref{secdiscrser} recalls facts concerning the
principal series of representations of~$\PSL_2(\RR)$.
In Section~\ref{secauthf} we give the relation between
automorphic forms and invariant hyperfunctions.

The subject of Section~\ref{secgeodec} is the geodesic
decomposition of hyperfunctions associated to
automorphic forms with polynomial growth. As a
preparation we discuss in Section~\ref{secfourexp} the
Fourier expansion at a cusp, which we put at~$\infty$.
The condition of polynomial growth is only imposed
at~$\infty$. At other (not $\G$-equivalent) cusps the
growth may be arbitrary.

The geodesic decomposition represents the invariant
hyperfunction as a finite sum. An infinite sum is
considered in Section~\ref{secadd}. This is related to
the transfer operator, discussed in
Section~\ref{sectrop}.

In Section~\ref{seccov} we reformulate our results on
the universal covering group of $\PSL_2(\RR)$, and
give a cohomological interpretation of
$\G$-decompositions. In Section~\ref{secrho} we return
to the period function~$\ps$.

\ntst{Thanks.}I thank E.P.\,van den Ban,
J.J.\,Duistermaat, J.B.\,Lewis and D.\,Zagier for
their interest, help, and useful discussions.

Many of the ideas in this paper are present in the work
of Lewis, or have been the subject of our discussions
during Lewis's visits to Utrecht. Zagier has brought
the work of Lewis to my attention, and has shown
interest in this approach. Van den Ban showed me the
 argument in~\ref{decompR}. Over the years Duistermaat
 has repeatedly told me that invariant boundary forms
 should give insight into automorphic forms.

\nwsect{Automorphic forms}\label{secautforms}
\ntst{Examples of modular forms.}In
Section~\ref{secintro} we have already seen cuspidal
Maass forms. For $\re s>1$ the {\sl Eisenstein
series}\/ is given by
\begin{equation}
\label{Eisseries}
G(s;z) = \frac{\Gf(s)}{\p^s}
\sideset{}'\sum_{p,q\in\ZZ} \frac{y^s}{|qz+p|^{2s}}.
\end{equation}
The prime denotes that $(p,q)=(0,0)$ is omitted. {}From
this one can derive the following Fourier expansion:
\begin{equation}
\label{EisFourier}
G(s;z)= 2\Ld(2s) y^s + 2\Ld(2s-1) y^{1-s} +
\sum_{n\neq0} \frac{2\s_{2s-1}(|n|)}{|n|^s}
W_{0,s-1/2}(4\p |n|y) e^{2\p inx},
\end{equation}
with $\Ld(u) = \p^{-u/2} \Gf\left(\frac u2\right)
\z(u)$, $\z$ the zeta function of Riemann, and the
divisor sum $\s_w(m) = \sum_{d\divides m} d^w$. The
Fourier expansion defines $G(s;z)$ for all $s\in\CC$
except $s=0,\,1$. We have $G(s;-1/z)=G(s;z)$.

Holomorphic modular cusp forms occur for even
``weights'' $2k=12$ and $2k\geq 16$. They have a
Fourier expansion of the form $h(z)=\sum_{n=1}^\infty
c_n e^{2\p inz}$ and satisfy $h(-1/z)=z^{2k}h(z)$.

These various types of modular forms can be unified by
working on the group $\PSL_2(\RR)$.

\ntst{Notations.}\label{Gnots}Put $G\isdef
\PSL_2(\RR)$. Elements of~$G$ are indicated by a
representative in $\SL_2(\RR)$. So $\matr abcd$ and
$\matr{-a}{-b}{-c}{-d}$ denote the same element
of~$G$. Notations: $k(\th) \isdef
\matr{\cos\th}{\sin\th}{-\sin\th}{\cos\th}$, and $p(z)
\isdef \matc{\sqrt y}{x/\sqrt y} 0 {1/\sqrt y}$ for
$z\in\bhv$, $x=\re z$, $y=\im z$.

Conjugation by $\matr{-1}001 \in \PGL_2(\RR)$ gives an
outer automorphism $\matr abcd \mapsto j\matr abcd
\isdef \matr a{-b}{-c}d$ of~$G$; it is an involution.

Elements of~$G$ act on the upper half plane~$\bhv $ by
fractional linear transformations: $z\mapsto \matr
abcd \cdot z \isdef \frac{az+b}{cz+d}$. This action is
the restriction of the action of~$G$ on the complex
projective line $\prc\supset \bhv$ defined by the same
formula.

The Lie algebra $\glie_r$ of~$G$ is generated by $\HH
\isdef\matr100{-1}$, $\XX \isdef\matr0100$ and
$\YY\isdef\matr0010$. By $\glie$ we denote its
complexification $\glie_r\otimes_\RR\CC$. A convenient
basis of~$\glie$ is $\WW$, $\EE^+$, $\EE^-$, with
$\WW=\XX-\YY$ and $\EE^\pm \isdef \HH \pm
i\left(\XX+\YY\right) \in\glie$. The Casimir operator
is $\om \isdef -\frac14 \EE^+\EE^- +
\frac14\WW^2-\frac i2 \WW$; it determines a
bi-invariant differential operator on~$G$.

$N\isdef\vzm{\matr1x01}{x\in\RR}$ is a unipotent
subgroup of~$G$, and $A\isdef\vzm{p(iy)}{y>0} $ a real
torus of dimension~$1$. $P \isdef NA$ is a parabolic
subgroup of~$G$. The group $K \isdef
\vzm{k(\th)}{\th\in\RR\mod\p\ZZ}$ is a maximal compact
subgroup of~$G$. As Haar a measure on~$K$ we use $dk =
\frac1\p \, d\th$, with $k=k(\th)$.

\ntst{Discrete subgroup.}\label{Gamprop}We consider a
cofinite discrete subgroup~$\G$ of~$G$ with at least
one cuspidal orbit. By conjugation we arrange that
$\infty$ is a cusp of~$\G$, and that
$p\left(i+1\right)=\matr1101 $ generates the
subgroup~$\G_\infty$ of elements of~$\G$ that
fix~$\infty$. Note that $\G$ is allowed to have more
than one $\G$-orbit of cusps.

The fundamental example in this paper is the {\sl
modular group} $\Gmod\isdef \PSL_2(\ZZ)$. Here the set
of cusps is $ \prq$; it consists of one $\Gmod$-orbit.
The elements $\matr1101$ and $\matr01{-1}0=k(\p/2)$
generate~$\Gmod$.

\ntst{Automorphic forms.}By an automorphic form we mean
a function $u:G \rightarrow\CC$ that satisfies
\begin{enumerate}
\item[i)] $u(\g g)=u(g)$ for all $\g\in\G$,
\item[ii)] $u(gk(\th)) = u(g)e^{ir\th}$ for all
$k(\th)\in K$, for some $r\in 2\ZZ$, the {\sl weight},
\item[iii)] $\om u=s\left(1-s\right) u$ for some
$s\in\CC$, the {\sl spectral parameter}.
\end{enumerate}
Note that there are no growth conditions. This
definition is insensitive to the change $s\mapsto 1-s$
in the spectral parameter.

\ntst{From upper half plane to group.}Let $u$ be a
cuspidal Maass form as in~\ref{mafo}, and put
$u_0(p(z)k(\th)) \isdef u(z)$. It is not difficult to
check that $u_0$ is an automorphic form for~$\Gmod$
with weight~$0$ and eigenvalue $s\left(1-s\right)$.
The same holds for the Eisenstein series. We use the
same notation for $z\mapsto E(s;z)$ and $p(z)k(\th))
\mapsto E(z;p(z)k(\th))\isdef E(s;z)$.

To a holomorphic cusp form~$H$ of weight~$2k$ we
associate the function $ h(p(z)k(\th))\isdef y^k H(z)
e^{2ik\th}$. This is an automorphic form of
weight~$2k$ with eigenvalue $k-k^2$.

\ntst{}Each automorphic form is determined by the
function $z\mapsto u(p(z))$ on~$\bhv$, and satisfies
an elliptic differential equation. So it is a real
analytic function.

The Lie algebra acts by differentiation on the right.
For an automorphic form~$u$ with weight~$r$ and
spectral parameter~$s$ we have $\WW u=ir u$, and
$\EE^\pm u$ is an automorphic form of weight $r\pm 2$,
with the same spectral parameter. $\EE^\mp \EE^\pm u$
is always a multiple of~$u$. If $u$ is an automorphic
form on~$G$ with weight~$2k$, then the function
$z\mapsto y^{-k}u(p(z))$ is holomorphic if and only if
$\EE^- u=0$.

Automorphic forms for the group~$\Gmod$ are called {\sl
modular forms}.

\ntst{Reflection.}We define the involution~$\spgl$ on
functions on~$G$ by $\spgl f : g \mapsto f(j(g))$. It
satisfies $\spgl \circ \WW= -\WW\circ\spgl$, and
$\spgl\circ \EE^\pm = \EE^\mp \circ \spgl$.

If the involution~$j$ leaves $\G$ invariant (as is the
case for~$\Gmod$), then $\spgl$ preserves
$\G$-invariance on the left, and maps automorphic
forms to automorphic forms with the same spectral
parameter and opposite weight. The corresponding
eigenspace decomposition in weight~$0$ gives the
decomposition of Maass forms in even and odd ones.

\nwsect{Hyperfunctions}\label{sechf}
We consider the sheaves of hyperfunctions on the real
line~$\RR$ and on the circle $T\isdef \RR\bmod\p\ZZ$.
For a point of view that works in higher dimension we
refer to, e.g.,~\cite{Sk}.

\ntst{Holomorphic and analytic functions.}Let $\Okr$
denote the sheaf of holomorphic functions on the
complex projective line~$\prc$.

A {\sl real analytic}\/ function on an open set
$U\subset\RR$ is the restriction of an element of
$\Okr(W)$, where $W\supset U$ is an open set in~$\CC$,
that may depend on the function. So the sheaf~$\Akr$
of real analytic functions on~$\RR$ is the restriction
$\Okr|_\RR$. In the sequel we use `analytic' as
abbreviation of `real analytic', and say `holomorphic'
when we mean `complex analytic'.

\ntst{Hyperfunctions on~$\RR$.}(See~\cite{Sk}, \S1.1--3
for proofs and further information.) Let $U\subset
\RR$ be open, and choose some open $W\subset\CC$ such
that $U\subset W$. {\sl Hyperfunctions on}\/~$U$ are
elements of $\Okr\left(W\setminus U\right)\bmod
\Okr(W)$. This does not depend on the choice of~$W$.
We denote the linear space of hyperfunctions on~$U$
by~$\Bkr(U)$. This defines the sheaf~$\Bkr$ of
hyperfunctions on~$\RR$. Intuitively, a hyperfunction
represented by $g\in\Okr\left(W\setminus U\right)$ is
the jump in~$g$ when we cross~$U$. A more fancy
definition of the sheaf of hyperfunctions is $\Bkr =
\Hkr^1_\RR(\CC,\Okr)$ (sheaf cohomology).

 Multiplication of representatives makes $\Bkr$ into an
 $\Akr$-module. We map $\Akr(U)$ into $\Bkr(U)$ by
 sending $g\in \Okr(V)$, with $V\subset \CC$ open,
 $V\supset U$, to the hyperfunction represented by
$\th\mapsto g(\th)$ on $V\cap\bhv$ and $\th\mapsto 0$
on~$\ohv\isdef \vzm{z\in\CC}{\im z<0}$.

\ntst{Support.}The support~$\supp(\al)$ of a
hyperfunction $\al\in\Bkr(U)$ is the smallest closed
subset $C\subset U$ such that the restriction of~$\al$
to $U\setminus C$ is zero. A representative $g\in
\Okr(W)$ of~$\al$ extends holomorphically to the
points of~$U\setminus\supp(\al)$.

\ntst{Parting.}\label{decompR}
The sheaf~$\Bkr$ is flasque. This means that the
restriction maps $\Bkr(V) \rightarrow\Bkr(U)$ are
surjective for all open $U\subset V \subset\RR$.

Any $\al\in \Bkr(I)$ can be broken up at each point
$a\in I$: We can write $\al = \al_++\al_-$ with
$\al_\pm\in\Bkr(I)$, $\supp(\al_-) \subset [a,\infty)
\cap I$, $\supp(\al_+) \subset (-\infty,a] \cap I$.
Indeed, consider $\bt\in \Bkr\left(I\setminus\vz
a\right)$ that restricts to~$\al$ on $I \cap
(a,\infty)$ and to~$0$ on $I\cap (-\infty,a)$. The
flasqueness implies that there is an element of
$\Bkr(I)$ restricting to~$\bt$ on $I\setminus\vz a$.
This element we take as $\al_+$, and $\al_-\isdef
\al-\al_+$.

We call the decomposition $\al=\al_++\al_-$ a {\sl
parting}\/ of~$\al$ at~$a$. It is well defined in the
stalk~$\Bkr_a$. Al partings of~$\al$ at~$a$ are
obtained by replacing $\al_\pm$ by $\al_\pm \pm\n$,
where $\n\in\Bkr(I)$ satisfies $\supp \n \subset \vz
a$.

\begin{figure}
\begin{center}
\scriptsize
\setlength\unitlength{.6cm}
\begin{picture}(12,3)(-6,-1.5)
\put(-6,0){\line(1,0){12}}
\put(0,.2){$\supp(\al)$}
\put(-3,1.3){$C$}
\thicklines
\qbezier(-3,0)(-3,1.5)(-2,1.5)
\qbezier(-3,0)(-3,-1.5)(-2,-1.5)
\qbezier(4,0)(4,1.5)(3,1.5)
\qbezier(4,-0)(4,-1.5)(3,-1.5)
\put(-2,1.5){\line(1,0){5}}
\put(-2,-1.5){\line(1,0){5}}
\put(-2,0){\line(1,0){5}}
\put(.65,-1.5){\vector(-1,0){.3}}
\put(.35,1.5){\vector(1,0){.3}}
\end{picture}
\end{center}
\caption[p]{Contour used to define the duality
in~\ref{dualR}} \label{figCR}
\end{figure}
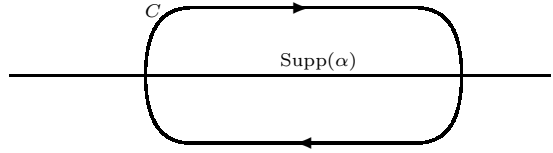

\ntst{Duality.}\label{dualR}Let $\Bkr_b(I)\isdef
\vzm{\al\in \Bkr(I)}{\supp(\al) \text{ is bounded}}$
be the space of hyperfunction on the open interval~$I$
with compact support. A duality between $\Akr(\RR)$
and $\Bkr_b(\RR)$ is given by
\[ \left\langle\ph, \al \right\rangle \isdef \int_C
\ph(\th) g(\th)\, \frac{d\th}\p \]
for $\ph \in \Akr(\RR)$, $g \in \Okr\left( W\setminus
\supp(\al)\right)$ a representative
of~$\al\in\Bkr_b(\RR)$ and $C$ any contour around
$\supp(\al)$ contained in~$W$ and in the domain of a
holomorphic function extending~$\ph$, see
Figure~\ref{figCR}.%
The use of the variable~$\th$ on~$\CC$, and the
measure~$\frac{d\th}\p$ will become clear in~\ref{T}.

\ntst{Reflection.}$\spgl \ph(\th)\isdef \ph(-\th)$
defines an involution~$\spgl$ in $\Akr(\RR)$. We
define the involution $\spgl$ in $\Bkr(\RR)$ by the
action $g\mapsto -\spgl g$ on representatives. In this
way $\spgl$ respects the injection $\Akr(\RR)
\rightarrow \Bkr(\RR)$ and satisfies $\left\langle
\spgl\ph, \spgl\al \right\rangle = \left\langle
\ph,\al \right\rangle$.
\medskip

\ntst{The circle~$T$.}The fact that $\Bkr$ is a sheaf
means that the definition of hyperfunctions is local,
and can be transferred to any real manifold of
dimension~$1$. We need hyperfunctions on the circle~$T
\isdef \RR\bmod\p\ZZ$. There are many ways to embed
$T$ into~$\prc$, for example, by $\th\mapsto
e^{2i\th}$ we view $T$ as the unit circle in~$T$. In
the sequel it is convenient to identify $T$ to the
common boundary of the upper half plane~$\bhv$ and the
lower half plane~$\ohv \isdef \vzm{z\in\CC}{\im z<0}$.
This we accomplish by the map $\pr:
\RR\rightarrow\prc: \th\mapsto \cot\th$.

We use the standard cyclic ordering on~$T$ induced by
$\RR\subset T$. Intervals in~$T$ are formed with
respect to this ordering: $[-1,1]$ is the same as the
corresponding interval in~$\RR$, but $[1,-1]=
[1,\infty) \cup \vz\infty \cup (-\infty,-1]$. The
map~$\pr$ is strictly decreasing. In~\ref{realindrpr}
we shall explain why we do not choose the increasing
map $\th\mapsto -\cot\th$.

\ntst{Analytic functions and hyperfunctions.} We define
the sheaves $\Akr_T$ of analytic functions on~$T$,
and~$\Bkr_T$ of hyperfunctions on~$T$ in the same way
as above: analytic functions on $U\subset T$ are the
restrictions of holomorphic functions on some open set
in~$\prc$ containing~$U$, and hyperfunctions on~$U$
are the elements of $\Okr\left(W\setminus
U\right)\bmod \Okr(W)$ for any fixed $W\supset U$.

\begin{figure}[t]
\setlength\unitlength{.7pt} \scriptsize
\hspace*{\fill}
\begin{picture}(170,170)(-80,-80)
\put(-85,0){\line(1,0){170}}
\put(90,-2){$T$}
\thicklines
\put(0,11){\line(1,0){75}}
\put(75,-11){\vector(-1,0){75}}
\put(-75,11){\vector(1,0){75}}
\put(0,-11){\line(-1,0){75}}
\put(75,-11){\line(0,-1){50}}
\put(-75,-11){\vector(0,-1){50}}
\put(-75,11){\line(0,1){50}}
\put(75,11){\vector(0,1){50}}
\qbezier(-75,61)(-75,80)(-40,80)
\qbezier(-75,-61)(-75,-80)(-40,-80)
\qbezier(75,61)(75,80)(40,80)
\qbezier(75,-61)(75,-80)(40,-80)
\put(-40,80){\line(1,0){80}}
\put(-40,-80){\line(1,0){80}}
\put(-50,18){$C_+$}
\put(-50,-22){$C_-$}
\end{picture}
\hfill
\begin{picture}(170,170)(-80,-80)
\qbezier(40,0)(40,40)(0,40)
\qbezier(0,40)(-40,40)(-40,0)
\qbezier(-40,0)(-40,-40)(0,-40)
\qbezier(0,-40)(40,-40)(40,0)
\put(43,-15){$T$}
\put(83,-15){$C_-$}
\put(5,-15){$C_+$}
\thicklines
\qbezier(80,0)(80,80)(0,80)
\qbezier(0,80)(-80,80)(-80,0)
\qbezier(-80,0)(-80,-80)(0,-80)
\qbezier(0,-80)(80,-80)(80,0)
\qbezier(25,0)(25,25)(0,25)
\qbezier(0,25)(-25,25)(-25,0)
\qbezier(-25,0)(-25,-25)(0,-25)
\qbezier(0,-25)(25,-25)(25,0)
\put(80,10){\vector(0,-1){20}}
\put(25,-5){\vector(0,1){10}}
\end{picture}
\hspace*{\fill} \caption[y]{The contours $C_+$ and
$C_-$ in the $\tau$-plane (left), and the $w$-plane
(right), with $w=\frac{\tau-i}{\tau+i}=e^{-2i\th}$. Of
course the contours at the left and the right are not
equal, but homotopic. }\label{Ccont}
\end{figure}
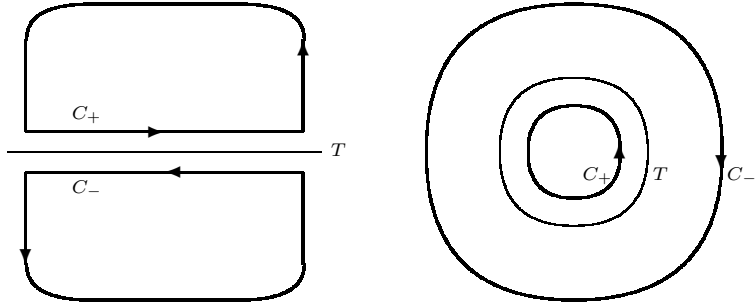

\ntst{Duality.}There is a duality between $\Akr_T(T)$
and $\Bkr_T(T)$ given by
\[ \langle f,\al\rangle \isdef \int_{C_+} f(\tau)
g(\tau) \frac{d\tau}{\p(1+\tau^2)} + \int_{C_-}
f(\tau) g(\tau) \frac{d\tau}{\p(1+\tau^2)},\]
for $\ph\in \Akr_T(T)$, $g$ a representative of $\al
\in \Bkr_T(T)$, and $C_\pm$ contours in the
intersections of the domains of representatives; see
Figure~\ref{Ccont}, and note the orientation.%

\ntst{Representative.}\label{repres}The hyperfunction
$\al\in \Bkr_T(T)$ can be recovered from the linear
form $\ph \mapsto \left\langle \ph, \al \right\rangle$
on~$\Akr_T(T)$. For each $\tau_0 \in \CC\setminus\RR$
the function $h_{\tau_0}:\tau \mapsto \frac {-i}2
\frac{1+\tau_0\tau}{\tau-\tau_0}$ determines an
element of~$\Akr_T(T)$. One can check that $g(\tau_0)
\isdef \left\langle h_{\tau_0}, \al \right\rangle$ is
the unique representative of~$\al$ that is holomorphic
on $\CC\setminus\RR$ and satisfies $g(i)+g(-i)=0$.

\ntst{Reflection.}The reflection on~$\RR$ corresponds
to the reflection induced by~$\tau\mapsto -\tau$
in~$\prc$. It defines an involution~$\spgl$ in
$\Akr_T(T)$ and $\Bkr_T(T)$ (again use an additional
minus sign in the representatives of hyperfunctions).
This involution respects the embedding
$\Akr_T\rightarrow\Bkr_T$ and leaves the duality
invariant.

\ntst{Basis.}For each $r\in\ZZ$ we define
$\ph_{2r}:\tau \mapsto \left( \frac{\tau+i}{\tau-i}
\right)^r$ in~$\Akr_T(T) \subset \Bkr_T(T)$. We have
$\left\langle \ph_{2r}, \ph_{2q} \right\rangle
=\dt_{q+r}$.

Any $\ph \in \Akr_T(T)$ has an expansion $\ph = \sum
c_r \ph_{2r}$, corresponding to a Laurent expansion
converging on an annulus of the form $p^{-1} < \left|
\frac{\tau+i}{\tau-i} \right| <p$ for some $p>1$.

Any $\al\in \Bkr_T(T)$ can be represented by a function
that is of the form $-\frac12 d_0 - \sum_{r=1}^\infty
d_r \left( \frac{\tau+i}{\tau-i} \right)^r$ on~$\ohv$,
and $\frac12 d_0 + \sum_{r=-1}^{-\infty} d_r \left(
\frac{\tau+i}{\tau-i} \right)^r$ on~$\bhv$. The
condition on the coefficients is
\begin{equation}
\label{contest} d_r = \oh(p^{|r|})\quad \text{ as
$|r|\rightarrow\infty$ for each $p>1$.}
\end{equation}

This shows that $\Bkr_T(T)$ corresponds to the space of
linear forms $L:\Akr_T(T) \rightarrow \CC$ that
satisfy $L\ph_{2r} = \oh\left( p^{|r|}\right)$ as
$|r|\rightarrow\infty$ for each $p>1$.\label{basest}

\nwsect{Principal series of
representations}\label{secdiscrser}
\ntst{Induced representation.}For each $\n\in\CC$ we
denote by $\p_\n$ the representation of~$G$ by right
translation in the space~$M^\n$ of classes of
functions $f:G\rightarrow\CC$ satisfying
$f(p(z)g)=y^{(1+\n)/2}$ and $\int_K
|f(k)|^2\,dk<\infty$. This gives the induced
representation of~$G$ corresponding to a character of
the parabolic subgroup~$P=\vz{\matc\ast\ast0\ast}$.

$M^\n$ is the Hilbert space $L^2(K,dk)$ with a
$G$-action~$\p_\n$ depending on~$\n$. This
representation is bounded. Under the pairing
$(f_1,f_2) \mapsto \int_K f_1(k) f_2(k)\,dk$ the
representations $(\p_\n,M^\n)$ and
$(\p_{-\n},M^{-\n})$ are dual to each other. See,
e.g., \cite{sl2r}, Chap.~III,~\S2.

Usually, the letter $H$ is used to indicate these
spaces. We employ~$M$ to avoid confusion with
cohomology groups.

\ntst{Realization of the induced
representation.}\label{realindrpr}The elements
of~$M^\n$ are sections of a line bundle over
$P\backslash G \cong K$. Here we view $P\backslash G$
as $T\subset\prc$, the boundary of~$\bhv$. We identify
$Pk(\th)$ with $\tau=\cot\th$. The right translation
in $P\backslash G$ by $g=\matr abcd\in G$ corresponds
to $\tau \mapsto g^{-1}\cdot\tau =
\frac{d\tau-b}{-c\tau+d}$.

If we would have chosen $Pk(\th) \mapsto -\cot\th$,
then the action would correspond to $\tau \mapsto
j(g)^{-1} \cdot\tau = \frac{d\tau+b}{c\tau+a}$. The
presence of~$j$ in this formula we dislike so much,
that we accept that $\th\mapsto \cot\th$ inverts the
order.

In terms of the variable~$\tau$ we find:
\begin{align*}
\p_\n\matr abcd \ph(\tau)
&= \left(\frac{1+\tau^2 }{ (c\tau-a)^2 + (d\tau-b)^2 }
\right)^{(1+\n)/2} \ph\left(
\frac{d\tau-b}{-c\tau+a}\right),
\\
\left\langle f_1, f_2 \right\rangle &= \frac1\p\int_T
f_1(\tau) f_2(\tau) \, \frac{d\tau}{1+\tau^2}.
\end{align*}

\ntst{Reflection.}The reflection~$\spgl$ in $M^\n$
considered as a space of functions on~$G$ corresponds
to the reflection $\spgl\ph(\tau)=\ph(-\tau)$ in the
functions on~$T$. It satisfies $\spgl \p_\n(g) =
\p_\n(j(g)) \spgl$.

\ntst{Analytic functions and hyperfunctions.}The
formula defining $\p_\n(g)$ in the functions on~$T$
preserves analyticity. Let $M_\om^\n\isdef \Akr_T(T)$
be the space of analytic functions on~$T$. So
$\left(\p_\n, M_\om^\n\right)$ is an algebraic
subrepresentation of~$(\p_\n, M^\n)$. The factor
$\left(\frac{1+\tau^2 }{ (c\tau-a)^2 + (d\tau-b)^2 }
\right)^{(1+\n)/2}$ is holomorphic on a neighborhood
of~$T$ in~$\prc$.

The formula defining $\p_\n(g)$ also makes sense when
applied to representatives of hyperfunctions. Let
$M^\n_{-\om} \isdef \Bkr_T(T)$ be the space of
hyperfunctions on~$T$. This gives a representation
$(\p_\n, M^\n_{-\om})$. We have $\left\langle
\p_\n(g)\ph,\al\right\rangle = \left\langle
\ph,\p_{-\n}(g)^{-1} \al \right\rangle$ for $\ph\in
M^\n_\om$ and $\al\in M^{-\n}_{-\om}$.

In between $M^\n_\om$ and $M^\n$ there is the
$\p_\n(G)$-invariant space $C^\infty(T)$. Its dual,
the space of distributions, sits between $M^\n$
and~$M^\n_{-\om}$. We do not consider these spaces in
this paper.

\ntst{$K$-finite vectors.}\label{diffrel}All elements
of~$M^\n_\om$ are differentiable vectors of $(\p_\n,
M^\n)$. The action of~$\glie$ satisfies $d\p_\n(\WW)
\ph_{2r} = 2ir\ph_{2r} $ and $d\p_\n(\EE^\pm)\ph_{2r}
= \left(1+\n\pm 2r \right) \allowbreak \ph_{2r\pm2} $.
The reflection satisfies $\spgl\ph_{2r}=\ph_{-2r}$.

Let $M^\n_K \subset M^\n_\om$ be the space of finite
linear combinations of the~$\ph_{2r}$. It is invariant
under $d\p_\n(\glie)$. The $(\glie,K)$-modules
$(d\p_\n,M^\n_K)$ have been classified, see, e.g.,
\cite{sl2r}, Chap.~VI, \S5. We note the following
facts:

$(d\p_\n, M^\n_K)$ is irreducible if and only if $\n\in
\CC\setminus\left(1+2\ZZ\right)$. In this case
$(\p_\n,M^\n_K)$ and $(\p_{-\n},M^{-\n}_K)$ are
isomorphic. The isomorphism is determined up to a
factor. We choose $\iota(\n):M^\n_K\rightarrow
M^{-\n}_K$ given by
\begin{equation}
\label{iotadef}
\iota(\n) \ph_{2r} = \left( \frac{1-\n}2\right)_{|r|}
\left( \frac{1+\n}2\right)_{|r|}^{-1}.
\end{equation}
We have $\iota(-\n)\iota(\n)=1$.

If $\n>0$ is odd, then $M^\n_K$ has two irreducible
subspaces $D_K^+\left(\n+1\right) \isdef
\bigoplus_{2r\geq \n+1} \CC\ph_{2r}$ and
$D_K^-\left(\n+1\right) \isdef \bigoplus_{2r\leq
-\n-1} \CC\ph_{2r}$. These are the discrete series
representations.

If $\n<0$ is odd, then there is the irreducible finite
dimensional subspace $E\left(\n+1\right) \isdef
\bigoplus_{|2r| < 1-\n} \CC\ph_{2r}$.

For $k\in\NN \isdef \ZZ_{\geq1}$ there are the
following exact sequences of $(\glie,K)$-modules:
\[ \widearray{ccccccccc} 0 & \rightarrow & D^+_K(2k)
\oplus D^-_K(2k)
&\rightarrow & M^{2k-1}_K
&\rightarrow & E(2k)
&\rightarrow & 0 \\
0 &\rightarrow & E(2k)
&\rightarrow & M_K^{1-2k}
&\rightarrow & D^+_K(2k) \oplus D^-_K(2k)
&\rightarrow
& 0
\end{array}
\]
The homomorphisms are unique up to a factor. A possible
choice for the first one is the continuation
of~$\iota(\n)$ to $\n=2k-1$, and for the second one
$\res_{\n=1-2k} \iota(\n)$.

\ntst{Extension.}Let
$\n\in\CC\setminus\left(1+2\ZZ\right)$. The factor
$\left( \frac{1-\n}2\right)_{|r|} \left(
\frac{1+\n}2\right)_{|r|}^{-1}$ has polynomial growth
as $|r|\rightarrow\infty$. This implies that we have
extensions $\iota(\n): M^\n_\om \rightarrow
M^{-\n}_\om$ and $\iota(\n): M^\n_{-\om} \rightarrow
M^{-\n}_{-\om}$. These isomorphisms respect the
$G$-action, and satisfy $\left\langle \ph,\iota(\n)
\al \right\rangle = \left\langle \iota(\n)\ph, \al
\right\rangle$.

\ntst{The space $E(2k)$.}Let $k\in\NN$. Multiplication
by $\left(\tau^2+1\right)^{k-1}$ gives a bijection
from $E\left(2k\right)$ onto the polynomials in~$\tau$
of degree at most $2k-2$. The action $\p_{1-2k}\matr
abcd$ corresponds to $F\mapsto F|_{2-2k} g^{-1}$,
where $\left( F|_{2r} \matr abcd \right)(\tau) \isdef
\left(c\tau+d\right)^{-2r}
F\left(\frac{a\tau+b}{c\tau+d} \right)$.
\label{polact}

\ntst{}\label{topol}For $k\in\NN$ we give a map $\al
\mapsto \al^{\langle 2k \rangle}$ from
$M^{2k-1}_{-\om}$ to the polynomials of degree at most
$2k-2$, extending the composition $M^{2k-1}_K
\rightarrow E(2k) \rightarrow \text{(polynomials)}$,
by
\[ \al^{\langle 2k \rangle}(X) \isdef \left\langle h_X,
\al \right\rangle, \text{ with } h_X(\tau) \isdef
-(2i)^{2k-2} \left( \tau^2+1\right)^{1-k}
(\tau-X)^{2k-2}. \]
A computation shows that $\al \mapsto
\al^{\langle2k\rangle}$ respects the $G$-action, and
vanishes on $\sum d_r \ph_{2r} \in M^{2k-1}_{-\om}$
with $d_r=0$ for $|r|\geq k$.

\nwsect{Automorphic hyperfunctions}\label{secauthf}

\ntst{Poisson integral.}In Theorem~3 of~\cite{He72},
Helgason shows that all eigenfunctions of the
Laplacian in the non-Euclidean plane can be described
as the Poisson integral of a hyperfunction. In our
notation this result states that each eigenfunction~$F
\in C^\infty(G/K)$ of the Casimir operator~$\om$ with
eigenvalue $\frac14\left(1-\n^2\right)$ can be be
written as
\[ F(p(z)k(\th)) = \left\langle \p_{-\n}(p(z))\ph_0,
\al\right\rangle\]
for some $\al \in M^\n_{-\om} = \Bkr_T(T)$.
In~\cite{Sk}, \S5, especially~5.5, we see that this is
only a very special case of general results for
symmetric spaces.
\ntst{Hyperfunctions and $(\glie,K)$-modules.}Let $\al
\in M^\n_{-\om}$. For each $\ph\in M^{-\n}_\om$ we put
$T_\al \ph(g) \isdef \left\langle \p_{-\n}(g)\ph, \al
\right\rangle$. This defines a linear map $T_\al :
M^{-\n}_\om \rightarrow C^\infty(G)$ that intertwines
$\p_{-\n}$ with the action of~$G$ by right
translation.

If we restrict $T_\al$ to $M^{-\n}_K$ we get a
$(\glie,K)$-module in $C^\infty(G)$ that is isomorphic
to a quotient of~$M^{-\n}_K$. Conversely, Helgason's
proof can easily be generalized to show that each such
$(\glie,K)$-module is described by an unique~$\al \in
M^\n_{-\om}$.

Under this correspondence the property
$\p_\n(\g)\al=\al$ for some $\g\in G$ is equivalent to
$(T_\al \ph)(\g g) = (T_\al \ph)(g)$ for all $\ph \in
H^{-\n}_K$. Actually, it suffices to let $\ph$ run
through the~$\ph_{2r}$.

\ntst{Definition.}Let $A^\n_{-\om}(\G)$ be the space of
$\al\in M^\n_{-\om}$ that satisfy $\p_\n(\g)\al=\al$
for all $\g\in\G$. The elements of~$A^\n_{-\om}(\G) $
we call {\sl automorphic hyperfunctions}\/ for~$\G$.

If $\al$ is an automorphic hyperfunction, then $T_\al
M^{-\n}_K$ is a $(\glie,K)$-module consisting of
linear combinations of automorphic forms, and all
$(\glie,K)$-modules isomorphic to a quotient of
$M^{-\n}_K$ in which the weight vectors are
automorphic forms arise in this way.

If $j(\G)=\G$, then $\spgl$ maps $A^\n_{-\om}(\G)$ into
itself. We have $\spgl(T_\al\ph) = T_{\spgl\al}
(\spgl\ph)$.

\ntst{Hyperfunction for holomorphic automorphic
forms.}\label{hypholo}As an example we consider a
holomorphic automorphic form~$H$ for~$\G$ with even
weight~$2k$. So $H\left( \frac{az+b}{cz+d}\right) =
\left(cz+d\right)^{2k} H(z)$ for $z\in\bhv$ and $\matr
abcd\in\G$. We do not impose any condition at the
cusps of~$\G$, so $k$ may be negative.

Define $\al \in M^{2k-1}_{-\om}$ to be the
hyperfunction represented by the function~$g$ equal
to~$0$ on~$\ohv$ and given by $g(\tau) = (-1)^k 4^{-k}
\left(1+\tau^2\right)^k H(\tau)$ for $\tau\in\bhv$.
The transformation behavior of~$H$ under~$\G$ implies
that $\al \in A^{2k-1}_{\om}$. A computation shows:
\[ \left\langle \p_{1-2k}(p(z)) \ph_{2r}, \al
\right\rangle =\frac{(-1)^k 4^{-k}}\p y^{1-k}
\int_{C_+} H(\tau) (\tau-\bar z)^{k+r-1}
(\tau-z)^{k-r-1}\,d\tau. \]
This vanishes if $r<k$, and yields $y^k H(z)$ for
$r=k$. So the automorphic form on~$G$ corresponding
to~$H$ is equal to $T_\al \ph_{2k}$. The element
$T_\al \ph_{2k}$ is a lowest weight vector in the
$(\glie,K)$-module it generates. If $2k\geq2$, this
$(\glie,K)$-module is isomorphic to $D^+(2k)$; it is
the image of $M^{1-2k}_K \rightarrow D^+_K(2k)\oplus
D^-_K(2k) \rightarrow D^+_K(2k)$. If $k\leq0$, the
$(\glie,K)$-module is not irreducible. It is
isomorphic to $M^{1-2k}_K \bmod
D^-_K\left(2-2k\right)$.

\ntst{Maass forms.}Any automorphic form of weight zero
generates a $(\glie,K)$-module that is the quotient of
some $M^\n_K$. Helgason's result quoted above shows
that these automorphic forms all arise from
automorphic hyperfunctions. In
Section~\ref{secfourexp} we shall give an explicit
construction of the hyperfunction corresponding to
automorphic forms with polynomial growth at the
cusp~$\infty$.

If the eigenvalue is $s\left(1-s\right)$ with
$s\not\in\ZZ$, then both $\n=2s-1$ and $\n=1-2s$ are
possible; the resulting automorphic hyperfunctions are
unique, and are related by~$\iota(\n)$. If $s\in\ZZ$,
only one of these choices will work, the hyperfunction
need not be unique.

If $j(\G)=\G$, and $\al \in A^\n_{-\om}(\G)$
corresponds to the Maass form~$u$, then $\spgl\al$
corresponds to the Maass form~$z\mapsto u(-\bar z)$.

\ntst{Eisenstein series in the domain of absolute
convergence.}\label{Eishf}For $\re s>1$ we define
\[ h_s(\tau) \isdef \frac{-i}2 \p^{-s} \Gf(s)
\sideset{}'\sum_{p,q\in\ZZ} \left( p^2+q^2\right)^{-s}
\frac{p\tau-q}{q\tau+p}.\]
This converges absolutely for all $\tau
\in\CC\setminus\RR$. The convergence is uniform on
compact sets in $\bhv\cup\ohv$. Let $\e_s^\ast$ be the
hyperfunction on~$T$ represented by~$h_s$. The
integral for $\left\langle \ph,\e^\ast_s\right\rangle$
can be evaluated term by term. For each term the
integrand has only one pole on~$T$, at $\tau=-\frac
pq$. We obtain $\left\langle
\ph,\e^\ast_s\right\rangle =\p^{-s} \Gf(s)
\sideset{}'\sum_{p,q} \ph\left(-\frac pq\right)
\left(p^2+q^2\right)^{-s}$. If we take $\ph(\tau) =
\p_{2s-1}(p(z)) \ph_0(\tau) = y^s \allowbreak \left(
\frac{\tau^2+1}{(\tau-z)(\tau-\bar z)} \right)^s$,
then we find $G(s;z)$. For $s\not\in \ZZ$ this
determines the hyperfunction uniquely. So for $\re
s>1$, $s\not\in\ZZ$, we have $\e^\ast_s \in
 A^{1-2s}_{-\om}(\Gmod)$. The relation giving the
 equivalence of $\p_{1-2s}(\g)h_s$ and $h_s$ for
 $s\not\in\ZZ$ extends to $s\in1+2\ZZ$, $s\geq2$. Hence
 $\e^\ast_s\in A^{1-2s}_{-\om}(\Gmod)$ for all~$s$ with
$\re s>1$. It corresponds to the Eisenstein series,
and satisfies $\spgl\e^\ast_s=\e^\ast_s$.

\ntst{Exponentially growing Poincar\'e series.}In their
construction and meromorphic continuation of
Poincar\'e series, Miatello and Wallach, \cite{MW},
explicitly give the linear form corresponding to an
automorphic hyperfunction. Their context is much wider
than ours: Lie groups with real rank one. Their
Poincar\'e series have in general exponential growth
at a cusp.

\ntst{Question.}Are there automorphic forms that
generate a $(\glie,K)$-module which is not the
quotient of some $M^{-\n}_K$?

In the sequel we consider automorphic hyperfunctions as
the principal objects.

\ntst{Modular case.}For the modular group, automorphic
hyperfunctions are closely related to functions
satisfying~\vgl{MFE}. This is the subject of the
remaining part of this section.
Theorem~\ref{authftopsi} is the main result. The
intermediate result Proposition~\ref{alphtof} is valid
for for all cofinite discrete groups~$\G$ with cusps.

\ntst{Definition.}Let
$\n\in\CC\setminus\left(1+2\ZZ\right)$. We define
$\Psimod(\n)$ to be the linear space of holomorphic
functions $\psi:\CC\setminus(-\infty,0]\rightarrow\CC$
that satisfy
\begin{gather}
\label{MFEnu}
\ps(\tau) = \ps(\tau+1) + (\tau+1)^{-\n-1}
\ps\left(\frac\tau{\tau+1}\right),
\\
\label{limsom}
\lim_{\im\tau\rightarrow\infty} \left( \ps(\tau) +
\tau^{-1-\n}\ps\left(\frac{-1}\tau\right) \right) +
\lim_{\im\tau\rightarrow -\infty} \left( \ps(\tau) +
\tau^{-1-\n}\ps\left(\frac{-1}\tau\right) \right) =0.
\end{gather}
The existence of both limits is part of
condition~\vgl{limsom}. Equation~\vgl{MFEnu} is
equation~\vgl{MFE}, with $2s$ replaced by $\n+1$.

\begin{stel}{Theorem}\label{authftopsi}For each
$\n\in\CC\setminus\left(1+2\ZZ\right)$ there is a
bijective linear map $A^\n_{-\om}(\Gmod)
\rightarrow\Psimod(\n):\al\mapsto\ps_\al$.
\end{stel}
\bw{Remarks} The {\sl proof}\/ is given in
\ref{Fdef}--\ref{alphtof}. The map is the composition
of the maps described in Lemma~\ref{FkrtoPsikr} and
Proposition~\ref{alphtof}. For $\n\in1+2\ZZ$ the map
$\al\mapsto\ps_\al$ from automorphic hyperfunctions to
holomorphic functions on~$\CC\setminus(-\infty,0]$
that satisfy~\vgl{MFEnu} is well defined, but we have
no bijectivity.

Equation~\vgl{MFEnu} is essential in the definition
of~$\Psimod(\n)$. The limit condition~\vgl{limsom} is
a normalization, needed to obtain injectivity.

\ntst{Definitions.}\label{Fdef}The space $\Psimod(\n)$
is contained in the space~$\Ps(\n)$ of holomorphic
functions on the smaller domain $\CC\setminus\RR$ that
satisfy \vgl{MFEnu} and~\vgl{limsom}.

For $\n\in\CC$, let $ \Fkr(\n)$ be the space of
holomorphic functions
$f:\CC\setminus\RR\rightarrow\CC$ that satisfy
\begin{align}
\label{period1}f(\tau) &= f(\tau+1),\\
\label{limfcond} f(\tau)
&= \oh(1)\text{ as $|\im\tau|\rightarrow\infty$, and }
f(i\infty) + f(-i\infty)=0.
\end{align}

If a $1$-periodic holomorphic function on~$\hv^\pm$ is
$\oh(1)$ as $\pm\im\tau \rightarrow\infty$, it has a
power series expansion in $e^{\pm2\p i\tau}$. Hence
$f(\pm i\infty)$ makes sense; it is the constant term
in the expansion.

\begin{stel}{Lemma}\label{FkrtoPsikr}Let $\n\in\CC$,
$\n\not\in 1+2\ZZ$. The relations
\begin{align}
\label{ftopsi} \ps(\tau)
&= f(\tau) -\tau^{-1-\n}f\left(\frac{-1}\tau\right)
\\
 \label{psitof}
f(\tau) &= \frac1{1+e^{\mp \p i\n } } \left( \ps(\tau)
+ \tau^{-1-\n}\ps\left( \frac{-1}\tau\right) \right)
\quad\text{ for $\tau\in\hv^\pm$}
\end{align}
define a bijective linear map $ \Fkr(\n) \rightarrow
\Ps(\n): f\mapsto \ps$.
\end{stel}
\bw{Remarks}J.~Lewis has shown me these
transformations. See~\ref{modgrcoh} for a
cohomological interpretation.

If $\n\in 1+2\ZZ$, then \vgl{ftopsi} defines a map from
$\Fkr(\n) $ to the holomorphic functions
on~$\CC\setminus\RR$ that satisfy~\vgl{MFEnu}.
\bw{Proof}For $f\in \Fkr(\n)$ define $\ps$
by~\vgl{ftopsi}. Then \vgl{psitof} turns out to give
back~$f$. The periodicity~\vgl{period1} for~$f$ is
equivalent to~\vgl{MFEnu} for~$\ps$, and \vgl{limsom}
is a direct reformulation of~\vgl{limfcond}.

\ntst{Definition.}Let
$\n\in\CC\setminus\left(1+2\ZZ\right)$. We define
$\Fkrmod(\n)$ to be the subspace of $\Fkr(\n)$
corresponding to $\Psimod(\n) \subset \Ps(\n)$ under
the map $f\mapsto\ps$ of Lemma~\ref{FkrtoPsikr}.

The $f\in \Fkrmod(\n)$ are characterized by the
property that $\tau\mapsto f(\tau)-\tau^{-1-\n}
f\left(\frac{-1}\tau\right)$ extends holomorphically
across $(0,\infty)$. This is equivalent to
$\tau\mapsto f(\tau) - (-\tau)^{-1-\n}
f\left(\frac{-1}\tau\right) =
-(-\tau)^{-1-\n}\ps\left(\frac{-1}\tau\right)$ having
a holomorphic extension across $(-\infty,0)$.

For $\n\in 1+2\ZZ$ we define $\Fkrmod(\n)$ as the space
of $f\in \Fkr(\n)$ for which $\tau \mapsto f(\tau) -
\tau^{-1-\n}f\left(\frac{-1}\tau\right)$ extends
holomorphically to~$\CC\setminus\vz0$.

\begin{stel}{Proposition}\label{alphtof}Let $\n\in\CC$,
and $\G\subset G$ as in~\ref{Gamprop}. There is an
injective linear map $ A^\n_{-\om}(\G) \rightarrow
\Fkr(\n):\al \mapsto f_\al$ such that $\tau \mapsto
\left(1+\tau^2\right)^{(1+\n)/2} f_\al(\tau)$
represents the restriction of the hyperfunction~$\al$
to the open subset $T\setminus\vz\infty$ of~$T$.

If $\G=\Gmod$, then the image of $ A^\n_{-\om}(\G)
\rightarrow \Fkr(\n): \al \mapsto f_\al$ is equal
to~$\Fkrmod(\n)$.
\end{stel}
\bw{Remarks}Here we do not need to exclude $\n\in
1+2\ZZ$. See~\ref{holofalph} for the case
corresponding to holomorphic automorphic forms.
\bw{Proof}The restriction~$\al_0$ of $\al$ to the open
$N$-orbit $T\setminus\vz\infty \subset T$ is
represented by a function~$g$ that is holomorphic on
at least the strips $0< |\im\tau|< \e$ for some
$\e>0$. Take $\e<1$. Then $F: \tau \mapsto
\left(1+\tau^2\right)^{-(1+\n)/2} g(\tau) $ is also
holomorphic on these strips. The invariance of~$\al_0$
under $\p_\n\matr 1101$ implies that
$F\left(\tau-1\right) = F(\tau) + q(\tau)$, with $q$
holomorphic on $|\im\tau|<\e$. So $F$ represents a
hyperfunction on~$\RR$ that is invariant under the
translations $\tau \mapsto \tau+n$ with $n\in\ZZ$. It
determines a hyperfunction on the circle, and that
hyperfunction can be represented by a function
holomorphic on the complement of the circle in~$\prc$.
This function is unique up to an additive constant. So
$F$ can be replaced by the unique function $f_\al$ of
the form
\begin{equation}
\label{falphdef}
f_\al(\tau) =
\begin{cases}
\frac12 A_0(\al) + \sum_{n=1}^\infty A_n(\al) e^{2\p
in\tau}
\quad & \text{ for $\tau\in\bhv$,}\\
-\frac12 A_0(\al) - \sum_{n=1}^\infty A_{-n}(\al)
e^{-2\p in\tau}
\quad & \text{ for $\tau\in\ohv$.}
\end{cases}
\end{equation}
The function $\tau \mapsto \left(
1+\tau^2\right)^{(1+\n)/2} f_\al(\tau)$ is holomorphic
on $\left( \bhv\setminus i[1,\infty) \right) \cup
\left( \ohv\setminus (-i)[1,\infty) \right)$. It
represents $\al_0$. If $f_\al$ would vanish, then
$\al_0=0$. As $\infty$ cannot be a fixed point of the
whole group~$\G$, this implies $\al=0$. Hence
$\al\mapsto f_\al$ is injective.

It is clear that $f_\al \in \Fkr(\n)$. To see that it
is in~$\Fkrmod(\n)$ if $\G=\Gmod$, we note that $\tau
\mapsto \left( 1+\tau^{-2} \right)^{(1+\n)/2}
f_\al\left( \frac{-1}\tau \right)$ represents the
restriction of~$\al = \p_\n\matr01{-1}0 \al$ to
$T\setminus\vz0$. The gluing conditions on
$(0,\infty)$ and $(-\infty,0)$ between the
representatives are just the existence of holomorphic
extensions across $(0,\infty)$ and $(-\infty,0)$ that
characterize $\Fkrmod(\n)$ inside~$\Fkr(\n)$.

Let $\G=\Gmod$, and start with $f\in \Fkrmod(\n)$.
Clearly, $\tau\mapsto \left(1+\tau^2\right)^{(1+\n)/2}
\allowbreak f(\tau)$ represents a hyperfunction
$\bt_0\in \Bkr_T\left( T\setminus\vz\infty\right)$,
that satisfies $\p_\n\matr1101 \bt_0=\bt_0$. Put
$\bt_\infty \isdef \p_\n\matr01{-1}0 \bt_0 \in
\Bkr_T\left(T\setminus \vz0 \right)$. A representative
of~$\bt_\infty$ is $\tau \mapsto \left(1+\tau^{-2}
\right)^{(1+\n)/2} f\left(\frac{-1}\tau\right)$. The
definition of~$\Fkrmod(\n)$ implies that $\bt_0$
and~$\bt_\infty$ coincide on $(0,\infty)$ and
$(-\infty,0)$. So there exists $\bt\in M^\n_{-\om}$
restricting to~$\bt_0$ on $T\setminus\vz\infty$ and
to~$\bt_\infty$ on~$T\setminus\vz0$. Clearly
$\p_\n\matr01{-1}0 \bt=\bt$. We have to check the
invariance under the other generator $\matr1101$ of
the modular group. The support of
$\p_\n\matr1101\bt-\bt$ is contained in~$\vz\infty$.
On $|\tau|>2$, $\tau\not\in\RR$, this hyperfunction is
represented by
\begin{multline*}
\left( \frac{ \tau^2+1}{1+(\tau-1)^2}
\right)^{(1+\n)/2}
\left(1+(\tau-1)^{-2}\right)^{(1+\n)/2}
f\left(\frac{-1}{\tau-1} \right)\\
\qquad\hbox{} - \left(1+\tau^{-2}\right)^{(1+\n)/2}
f\left( \frac{-1}\tau\right)
\\
= \left(1+\tau^{-2}\right)^{(1+\n)/2} \left( \left(
\frac\tau{\tau-1} \right)^{1+\n} f\left(
\frac{-1}{\tau-1} -1 \right) - f\left( \frac{-1}\tau
+1 \right) \right)
\\
= - \left(1+\tau^{-2}\right)^{(1+\n)/2} \left( f\left(
{\textstyle\frac{\tau-1}\tau}\right) - \left(
{\textstyle\frac{\tau-1}\tau}\right)^{-1-\n} f\left(
\frac{-1}{ {\textstyle\frac{\tau-1}\tau}} \right)
\right).
\end{multline*}
The fact that $f\in \Fkrmod(\n)$ implies that the
quantity between brackets is holomorphic on a
neighborhood of ${\textstyle\frac{\tau-1}\tau} =1$.
This shows that $\p_\n\matr1101\bt-\bt$ vanishes on a
neighborhood of~$\infty$.

\ntst{Holomorphic automorphic
forms.}\label{holofalph}Let $H$ be a holomorphic
automorphic form~$H$ for~$\G$ of even weight~$2k$. It
has a Fourier expansion $H(z) =
\sum_{n=-\infty}^\infty a_n e^{2\p inz}$ converging
for $z\in\bhv$. If it is bounded at the cusp~$\infty$,
then the sum is over $n\geq0$. Meromorphy at the cusp
corresponds to a sum over $n\geq -N$ for some
$N\geq1$.

Suppose that $\al \in A^{2k-1}_{-\om}(\G)$ corresponds
to~$H$ as indicated in~\ref{hypholo}. {}From the
representative~$g$ in~\ref{hypholo} we subtract the
function~$p:\tau\mapsto (-1)^k4^{-k}
\left(1+\tau^2\right)^k \allowbreak P(\tau)$, with
$P(z)\isdef \frac12 a_0 + \sum_{n=-\infty}^{-1} a_n
e^{2\p inz}$. The function~$P$ is holomorphic
on~$\CC$, hence $g-p$ represents the restriction
of~$\al$ to $T\setminus\vz\infty$. We see that
\[ f_\al(\tau) = (-1)^k 4^{-k} \cdot\begin{cases}
\frac12 a_0 + \sum_{n=1}^\infty a_n e^{2\p
in\tau}\quad&\text{ for $\tau\in\bhv$,}\\
-\frac12a_0 - \sum_{n=1}^\infty a_{-n} e^{-2\p in\tau}
\quad&\text{ for $\tau\in\ohv$.}
\end{cases}
\]

\ntst{Other examples.}In~\ref{fcffs} we shall see that
for other automorphic hyperfunctions~$\al$ as well the
function~$f_\al$ is closely related to the Fourier
expansion of the automorphic forms associated
to~$\al$.%
\begin{table}
\[ \widearray{|c|c|c|}
\hline
&&\text{$F_nu_{2r}$ is a linear combination of:}\\
\hline
n=0 & \n\neq 0 & y^{(1+\n)/2}\text{ and }y^{(1-\n)/2}
\\
\cline{2-3}
&\n=0 & y^{1/2}\text{ and } y^{1/2}\log y \\
\hline
\multicolumn{2}{|l|}{n>0}
& W_{r,\n/2}(4\p ny)\\
\hline
\multicolumn{2}{|l|}{n<0}
& W_{-r,\n/2}(4\p|n|y)
\\
\hline
\end{array}
\]
\caption[p]{Basis elements for the spaces of Fourier
terms of automorphic forms with polynomial growth.}
\label{Fourierterms}
\end{table}

\nwsect{Fourier expansion}\label{secfourexp}In this
section we give an explicit description of a
representative of an automorphic hyperfunction that
has polynomial growth at the cusp~$\infty$. The main
result is Lemma~\ref{cuspdecomp}.

Some parts of this section are technical. Reading
\ref{systaf}--\ref{polgrwth}, the notations
in~\ref{falph0c}, and Lemma~\ref{cuspdecomp} suffices,
if one is willing to accept later on some results on
representatives of hyperfunctions.

\ntst{System of automorphic forms.}\label{systaf}We
work with an automorphic hyperfunction $\al\in
A^\n_{-\om}(\G)$. To it corresponds a system $(
u_{2r})_{r\in\ZZ}$ of automorphic forms given by
$u_{2r}(g) = \left\langle
\p_\n(g)\ph_{2r},\al\right\rangle$. This system
satisfies the differential equations
\begin{equation}
\label{differentrels}
\EE^\pm u_{2r} = (1 - \n \pm 2r)u_{2r\pm2}
\qquad \text{for $r\in\ZZ$.}
\end{equation}

\ntst{Polynomial growth.}\label{polgrwth}We say that an
automorphic hyperfunction~$\al$ has polynomial growth
at~$\infty$ if $u_{2r}(p(z)) =\oh(y^a)$ as
$y\rightarrow\infty$, uniformly in~$x$, for each
weight $2r\in2\ZZ$.

The hyperfunctions associated to Eisenstein series and
to cuspidal Maass forms have polynomial growth. The
Poincar\'e series studied in~\cite{MW} have in general
exponential growth at the cusp at which they are
defined, but not at $\G$-inequivalent cusps.

\ntst{Fourier expansion.}The cuspidal Maass forms and
Eisenstein series discussed in Sections~\ref{secintro}
and~\ref{secautforms} have been given by their Fourier
expansion at~$\infty$. In general we have:
\begin{align*}
u_{2r}(p(z)) &= \sum_{n\in\ZZ} e^{2\p inx} F_n
u_{2r}(y),\\
F_nu_{2r}(y) &\isdef \int_{x\in\RR\bmod\ZZ} e^{-2\p
inx} u_{2r}(p\left(iy+x\right) )\,dx.
\end{align*}
The {\sl Fourier terms}\/~$F_nu_{2r}$ satisfy a
differential equation with a two-dimensional solution
space; see, e.g., \cite{Maa}, Chap.~IV, \S2, or
\cite{fa}, \S4.1, 4.2, 4.4. For $n\neq0$ the condition
of polynomial growth imposes an additional condition.
See Table~\ref{Fourierterms} for the possibilities.

\ntst{Systems of Fourier terms.}The differential
equations~\vgl{differentrels} imply the same equations
for each of the systems~$(F_nu_{2r})_r$ separately. So
these systems are linear combinations of systems of
solutions of the corresponding differential equation.
We restrict ourselves to the case of polynomial
growth, so we are lead to systems of the functions in
Table~\ref{Fourierterms} that
satisfy~\vgl{differentrels}. In
Table~\ref{Fouriersystems} we give all
possibilities.%
\begin{table}
\[ \widearray{|l|l|l|r|}\hline
 n=0 & \n\not\in1+2\ZZ_{\geq0},\, \n\neq0&
 y^{(1-\n)/2},&{\bf a}\\
&& \Gf(\frac{1-\n}2+|r|)^{-1} \left(
\frac{1+\n}2\right)_{|r|} y^{(1+\n)/2}
&{\bf c}
\\
\cline{2-4}
&\n\in 1+2\ZZ,\, \n>0 & y^{(1-\n)/2},&{\bf a}
\\
&& \Gf(\frac{1-\n}2+|r|)^{-1} \left(
\frac{1+\n}2\right)_{|r|} y^{(1+\n)/2} ,
&{\bf c}
\\
&& \sign(r) \Gf(\frac{1-\n}2+|r|)^{-1} \left(
\frac{1+\n}2\right)_{|r|} y^{(1+\n)/2} & {\bf d}
\\
\cline{2-4}
&\n=0 & y^{1/2},&{\bf a}
\\
&& y^{1/2}\left(\log y + l_r \right)& {\bf e}\\
\hline
\multicolumn{2}{|l|}{n\neq0, \;\;\e=\sign n}&
\frac{(-1)^r}{ \Gf(\frac{1-\n}2+\e r)\vphantom{A_A} }
W_{\e r,\n/2}(4\p|n|y) e^{2\p inx}&{\bf b}
\\
\hline
\end{array}
\]
\caption[y]{Bases of the spaces of systems of Fourier
terms with polynomial growth.

For each system $2r\in 2\ZZ$ denotes the weight. We
give the function on~$G$ in the point~$p(z)$ with $z =
x+iy \in\bhv$. }\label{Fouriersystems}
\end{table}
For $n\neq0$ the space of Fourier terms with polynomial
growth has dimension~$1$. For $n=0$ the dimension
equals~$2$, except in the case of odd negative values
of~$\n$. Further remarks on the systems in the table:
\begin{enumerate}
\item[{\bf a.}] This is the standard form
of~$M^{-\n}_K$ as the induced representation from~$P$
to~$G$.
\item[{\bf b.}] If $\n$ is odd, some of the functions
in the system of Fourier term vanish. This is due to
the fact that the quickly decreasing Whittaker
functions span a $(\glie,K)$-module that is not a
quotient of either $M^\n_K$ or $M^{-\n}_K$ if
$\n\in1+\nobreak2\ZZ$.
\item[{\bf c.}] For negative, odd~$\n$ the functions
vanish if $|2r|\geq 1-\n$, due to the factor $\left(
\frac{1+\n}2\right)_{|r|}$. For odd, positive~$\n$ the
Gamma factor produces zeros for $|2r|\leq \n-1$.
\item[{\bf d.}] For odd positive~$\n$ the
equations~\vgl{differentrels} give no coupling between
$2r\geq \n+1$ and $2r\leq-1-\n$. This explains the
presence of three linearly independent systems.
\item[{\bf e.}] At $\n=0$ the systems {\bf a} and~{\bf
c} are linearly dependent. The $l_r$ are determined up
to an additive constant.
\end{enumerate}

\ntst{$N$-equivariant hyperfunctions.}In
\ref{muinf}--\ref{lambdadef} we describe each system
in Table~\ref{Fouriersystems} by means of a
hyperfunction that transform according
to~$\p_\n\matr1x01\bt = e^{-2\p inx}\bt$. Actually, in
Table~\ref{Fouriersystems} we have chosen factors that
seem natural when working with the functions
inTable~\ref{Fourierterms}. Other factors turn out to
be more natural when working with hyperfunctions. We
have summarized the results in the
Tables~\ref{Neqvarhf} and~\ref{supprestr}.

\ntst{Support in~$\infty$.}\label{muinf}Any
$\p_\n(N)$-invariant hyperfunction with
support~$\vz\infty$ has a representative $g$
holomorphic on a neighborhood of~$\infty$, with
$\infty$ itself deleted. The difference
$\p_\n\matr1x01 g -g$ has to have a holomorphic
extension to~$\tau=\infty$ for each real~$x$. An
analysis of the unique representative of the form
$g(\tau) = \left(1+\tau^{-2}\right)^{(1+\n)/2}
\sum_{k=1}^\infty r_k \tau^k$ leads to a solution that
is valid for all $\n \in\CC$:
\begin{equation}
\label{mudef}
\m\in M^\n_{-\om},\quad\text{represented by
}\tau\mapsto \frac{-i}2 \tau.
\end{equation}
The corresponding linear form on~$M^{-\n}_\om$ is $\ph
\mapsto \ph(\infty)$. This gives the system~{\bf a} in
Table~\ref{Fouriersystems}.

In general, all solutions are multiples of~$\m$. But if
$\n\in \ZZ_{\geq1}$, $\n\geq1$, the solution space has
dimension~$2$, and is spanned by~$\m$ and
\begin{equation}
\label{mustardef}
\m^\ast_\n,\quad \text{represented by }\tau \mapsto
\frac{-i}2 \tau^{\n+1}\left(1+\tau^{-2}
\right)^{(1+\n)/2}.
\end{equation}

A computation shows that $\left\langle
\p_{-\n}(p(z))\ph_{2r},\m^\ast_\n \right\rangle $ is
equal to
\begin{multline*}
\lim_{R\rightarrow\infty} \frac1\p \int_{|\tau|=R}
\frac{-i}2\tau^{\n+1}\left(1+\tau^{-2}\right)^{(1+\n)/2}
y^{(1-\n)/2}\\
\qquad\hbox{} \cdot \left(
\frac{1+\tau^2}{(\tau-z)(\tau-\bar
z)}\right)^{(1-\n)/2} \left( \frac{\tau-\bar
z}{\tau-z} \right)^r \frac{d\tau}{1+\tau^2}\\
= y^{(1-\n)/2} \cdot \text{ coefficient of $u^\n$ in }
\left(1-uz\right)^{(\n-1)/2-r} \left( 1-u\bar
z\right)^{(\n-1)/2+r} \displaybreak[0] \\
= y^{(1-\n)/2} \cdot (2iy\sign r)^\n
\frac{\Gf(\frac{\n+1}2+|r|)}{\Gf(\n+1)
\Gf(\frac{1-\n}2+|r|)} \displaybreak[0]\\
= \sqrt\p i^\n (\sign r)^\n \frac{ \left(
\frac{1+\n}2\right)_{|r|}} { \Gf(1+\frac\n2) \Gf(
\frac{1-\n}2 +|r|) } \, y^{(1+\n)/2}.
\end{multline*}
So $\m^\ast_\n$ determines a multiple of the
system~{\bf c} if $\n$ is even, and of system~{\bf d}
if $\n$ is odd.

\ntst{Integral over the open
$N$-orbit.}\label{intNorbit}Let $n\in\ZZ$. Put
$n=\pm|n|$, with $\pm=+$ if $n=0$. We define
$t_0(\n)\isdef \sqrt\p\Gf\left(-\n/2\right)^{-1}$, and
$t_n(\n)\isdef 1$ if $n\neq0$. If $\re\n<0$ the
following integral converges for each $\ph\in
\Akr_T(T)$, $n\in\ZZ$:
\begin{equation}\label{Nintegral}
\int_{-\infty}^\infty \ph(\tau) t_n(\n) e^{2\p
in\tau}\left(1+\tau^2\right)^{(1+\n)/2}
\frac{d\tau}{\p(1+\tau^2)}.
\end{equation}
For each $n\in\ZZ$ this defines a linear form
on~$\Akr_T(T)$, with value bounded by the supremum
norm of~$\ph$ on~$T$. So its values on the basis
elements~$\ph_{2r}$ are estimated by~$\oh(1)$.
In~\ref{basest} we see that the linear form is given
by a hyperfunction, which we call $\k_n(\n)$.

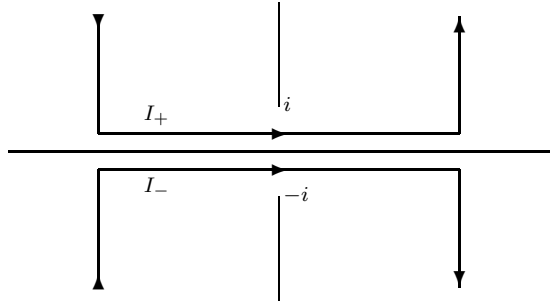
\begin{figure}[t]
\begin{center}
\footnotesize
\setlength\unitlength{.6cm}
\begin{picture}(10,5)(-5,-2.5)
\put(-6,0){\line(1,0){12}}
\put(0,1){\line(0,1){2.3}}
\put(0,-1){\line(0,-1){2.3}}
\put(0.1,.9){$i$}
\put(0.1,-1.1){$-i$}
\put(-3,.7){$I_+$}
\put(-3,-.9){$I_-$}
\thicklines
\put(-4,.4){\vector(1,0){4.2}}
\put(0,.4){\line(1,0){4}}
\put(4,.4){\vector(0,1){2.6}}
\put(-4,3){\line(0,-1){2.6}}
\put(-4,3){\vector(0,-1){.3}}
\put(-4,-.4){\vector(1,0){4.2}}
\put(0,-.4){\line(1,0){4}}
\put(4,-.4){\vector(0,-1){2.6}}
\put(-4,-3){\vector(0,1){.3}}
\put(-4,-3){\line(0,1){2.6}}
\end{picture}
\end{center}
\caption[u]{The contours $I_+$ and $I_-$.}\label{Icont}
\end{figure}

If $\pm n>0$ we can deform the path of integration
into~$I_\pm$, indicated in Figure~\ref{Icont}.%
The contour should be contained in the domain of~$\ph$.
The resulting integral converges for all $\n\in\CC$.
If we take the contour inside the region $c^{-1}<
\left| \frac{\tau+i}{\tau-i} \right|<c$ for some
$c>1$, we find that the value on~$\ph_{2r}$ is
$\oh(c^{|r|})$. So this integral extends the
definition of~$\k_n(\n)$ for $n\neq0$.

In all cases $\n\mapsto \left\langle \ph,
\k_n(\n)\right\rangle$ is holomorphic on the domain of
definition of $\k_n(\n)$. If we compute $\left\langle
\p_{-\n}(p(z))\ph_{2r}, \k_n(\n)\right\rangle$ by
means of~\vgl{Nintegral}, we find
\begin{alignat*}2
&\frac{(-1)^r \Gf(\frac{1-\n}2)}{ \Gf(\frac{1-\n}2+r)
\Gf(\frac{1-\n}2-r) } y^{(1+\n)/2}&\quad&\text{if
$n=0$,}\\
& \frac{(-1)^r (\p|n|)^{-(1+\n)/2} }{
\Gf(\frac{1-\n}2\pm r)} e^{2\p inx} W_{\pm r,\n/2}(4\p
|n|y)&\quad&\text{if $\pm n>0$.}
\end{alignat*}
Thus we have a multiple of system~{\bf c}, respectively
system~{\bf b}. This also shows that $\k_n(\n)$
behaves under $\p_\n\matr1x01$ according to the
character $\matr 1x01\mapsto e^{-2\p inx}$ of~$N$.

\ntst{Representative.}\label{kapnrepr}Let $\pm n>0$. As
we have seen in~\ref{repres}, we obtain a
representative of~$\k_n(\n)$ by
\[ g(\tau_0) = \frac1{2\p i} \int_{I_\pm}
\frac{1+\tau\tau_0}{\tau-\tau_0} e^{2\p in\tau}
\left(1+\tau^2\right)^{(\n-1)/2}\,d\tau. \]
In this integral $\tau_0$ is either in~$\hv^\mp$, or
inside the contour~$I_\pm$. We conclude that
$g(\tau_0)=\oh(1)$ as $|\im\tau_0|\rightarrow\infty$
uniformly for $\re\tau_0$ in compact sets.

In the case that $\tau_0\in\hv^\pm$, and
$|\im\tau_0|<1$, we deform the contour in such a way
that $\tau_0$ is outside it. This gives
\[ g(\tau_0) = \pm e^{2\p in\tau_0} \left(1+\tau_0^2
\right)^{(1+\n)/2} + \frac1{2\p i} \int_{I_\pm}
\frac{1+\tau\tau_0}{\tau-\tau_0} e^{2\p in\tau}
\left(1+\tau^2\right)^{(\n-1)/2}\,d\tau. \]

Thus we see that the restriction of~$\k_n(\n)$ to
$T\setminus\vz\infty$ can be represented by~$0$
on~$\hv^\mp$ and $\tau\mapsto \pm e^{2\p in\tau}
\left(1+\tau^2\right)^{(1+\n)/2}$ on $\hv^\pm
\setminus (\pm i)[1,\infty)$.

For $n=0$ and $\re \n < -\frac12$ we can draw similar
conclusions.%
\begin{figure}
\setlength\unitlength{.6cm}
\begin{center}
\scriptsize
\begin{picture}(10,6)(-5,-3)
\put(-5,0){\line(1,0){10}}
\put(0,1){\circle*{.05}}
\put(0,-1){\circle*{.05}}
\put(.1,.9){$i$}
\put(.1,-1.1){$-i$}
\put(1.2,3.3){$B_+$}
\put(-3,-2.1){$B_-$}
\put(-3,-.4){$A_+$}
\thicklines
\put(-4,0){\line(1,0){8}}
\qbezier(-4,0)(-4,3)(0,3)
\qbezier(0,3)(4,3)(4,0)
\qbezier(-4,0)(-4,-3)(0,-3)
\qbezier(0,-3)(4,-3)(4,0)
\put(.15,3){\vector(-1,0){.3}}
\put(-.15,-3){\vector(1,0){.3}}
\put(-.15,0){\vector(1,0){.3}}
\end{picture}
\end{center}
\caption[p]{Contours for the computation of
$\left\langle\ph, \k_0(\n) \right\rangle$
in~\vgl{kap0ABint}.}\label{AB}
\end{figure}

\ntst{Continuation of $\k_0(\n)$.}\label{kap0cont}Let
$\n\in\CC\setminus (-2\NN)$, and take
\begin{align*}
p(\tau)&\isdef
\begin{cases}
t_0(\n) \left(1+\tau^2\right)^{(1+\n)/2}&\text{ if
$\im\tau\geq0$, $\tau\not\in i[1,\infty)$,}\\
0&\text{ if $\tau\in\ohv$}
\end{cases}\\
q(\tau)&\isdef
\tau^{1+\n}\left(1+\tau^{-2}\right)^{(1+\n)/2} \frac{
\Gf(1+\frac\n2) e^{-\p i\n\sign(\im\tau)/2 } }{
2i\sqrt\p }\quad\text{for $\tau\not\in\RR\cup
i[-1,1]$.}
\end{align*}
(In~\ref{intNorbit} we have defined $t_0(\n)=\sqrt\p
\Gf(-\n/2)^{-1}$.)

 The functions $p$ and~$q$ define hyperfunctions
on~$T\setminus\vz\infty$, respectively $T\setminus0$.
They have been chosen in such a way that their
difference is holomorphic on $(0,\infty)$ and
on~$(-\infty,0)$. So together they define a
hyperfunction on~$T$, which we call $\bt$ for the
moment. We have the following integral representation
for each $\ph\in \Akr_T(T)$:
\begin{equation}
\label{kap0ABint}
\langle\ph, \bt \rangle = \int_{A_+} \ph(\tau)p(\tau)
\frac{d\tau}{\p(1+\tau^2)} + \sum_\pm \int_{B_\pm}
\ph(\tau) q(\tau) \frac{d\tau}{\p(1+\tau^2)},
\end{equation}
with contours as indicated in Figure~\ref{AB}. This
shows that $\langle\ph, \bt \rangle $ is holomorphic
in $\n\in\CC
\setminus(-2\NN)$.%
\begin{table}
\[ \widearray{|c|c|c|c|c|c|}
\hline
n&\n\in&\text{name}&
\text{representative}&\text{syst.}& \text{factor}
\\
\hline
0& \CC&\m& \tau\mapsto \frac{-i}2\tau & {\bf a}& 1
\\
\cline{2-6}
& \CC & \k_0(\n) & \text{see \ref{kapnrepr} and
\ref{kap0cont}} & {\bf c} & 1 \\
\cline{2-6}
& 1+2\ZZ_{\geq0}& \m^\ast_\n& \text{see
\vgl{mustardef}}
& {\bf d}&
\frac{i(-1)^{(\n-1)/2}\sqrt\p}{\Gf(1+\frac\n2)
\vphantom{A_{A_A}} }
\\
\cline{2-6}
& \vz0& \ld(0) & \text{see \ref{lambdadef}} & {\bf e} &
\frac2\p
\\
\hline
\neq0 & \CC& \k_n(\n) & \text{see \ref{kapnrepr}} &
{\bf b}& \left( \p|n|\right)^{-(1+\n)/2}
\\
\hline
\end{array}
\]
\caption[o]{$N$-equivariant hyperfunctions, satisfying
$\p_\n\matr 1x01\bt=e^{2\p inx}\bt$.

Each hyperfunctions determines a multiple of a system
of eigenfunctions in Table~\ref{Fouriersystems}. }
\label{Neqvarhf}
\end{table}

For $\re\n<0$ we can move off the contours~$B_\pm$ to
infinity, and obtain the integral in~\vgl{Nintegral}
with $n=0$. Thus we have extended the definition of
$\k_0(\n)$ to~$\CC$.

For $\n\in2\ZZ_{\geq0}$ the function~$p$ vanishes, and
the support of~$\k_0(\n)$ is~$\vz\infty$. In fact
$\k_0(\n) = (-1)^{\n/2} \p^{-1/2}
\Gf\left(1+\frac\n2\right) \m^\ast_\n$ if $\n\in2\NN$,
and $\k_0(0) = \p^{-1/2}\m$ (see~\ref{muinf}).

The $\p_\n(N)$-invariance, and the expression for
$\left\langle \p_{-\n}(p(z))
\ph_{2r},\k_0(\n)\right\rangle$ stay valid by
holomorphy.

\ntst{Logarithmic case.}\label{lambdadef}Let $\n=0$.
The function $\tau\mapsto \left( \frac {-i}\p\log\tau
-\frac{\sign\im\tau}2\right)\allowbreak \tau
\allowbreak \left(1+\tau^{-2} \right)^{1/2}$ on
$\CC\setminus\left(\RR\cup i[-1,1]\right)$ defines a
hyperfunction on a neighborhood of~$\infty$ that fits
nicely with the hyperfunction on~$T\setminus\vz\infty$
represented by~$\tau\mapsto
-\left(1+\tau^2\right)^{1/2}$ on~$\bhv\setminus
i[1,\infty)$ and $\tau\mapsto 0$ on~$\ohv$. We call
it~$\ld(0)$. It satisfies $\p_0
\matr1x01\ld(0)=\ld(0)$ for all $x\in\RR$. It is the
continuation to~$\n=0$ of the family $\ld:\n \mapsto
\frac2{\n\sqrt\p} \left( \k_0(\n)-\p^{-1/2}\m\right)$.

A computation shows that $\left\langle
\p_0(p(z))\ph_{2r},\ld(0) \right\rangle$ is equal to
$\frac2\p \sqrt y\log y$ plus a well defined but
complicated multiple of~$\sqrt y$.%
\begin{table}
\[ \widearray{|c|c|cl|}\hline
&\text{support}& p&\\
\hline
\m&\vz\infty& 0 &
\\
\hline
\k_n(\n)& T & e^{2\p in\tau}
\left(1+\tau^2\right)^{(1+\n)/2}& n\neq0\\
\hline
\k_0(\n)&T& \sqrt\p \Gf(-\n/2)^{-1}
\left(1+\tau^2\right)^{(1+\n)/2}&\n\not\in2\ZZ_{\geq0}
\\
&\vz\infty& 0& \n\in2\ZZ_{\geq0}
\\
\hline
\ld(0) & T & -\left(1+\tau^2\right)^{1/2}&
\\
\hline
\m^\ast_\n & \vz\infty& 0&
\\
\hline
\end{array}
\]
\caption[o]{Support and restriction to $T_0\isdef
T\setminus\vz\infty$ of the hyperfunctions in
Table~\ref{Neqvarhf}.

The restriction to~$T_0$ is represented by a function
$\tau\mapsto \pm p(\tau)$ on~$\hv^\pm\setminus\nobreak
(\pm i)[1,\infty)$, $\tau\mapsto0$ on~$\hv^\mp$. For
$\k_n(\n)$ the convention is $\pm n\geq0$, $\pm = +$
if $n=0$. } \label{supprestr}
\end{table}
\medskip

\ntst{Fourier terms of~$\al$.}\label{Ftal}Let $n\in\ZZ$
and consider a representative~$g$ of the automorphic
hyperfunction~$\al$. Define
\begin{align*} g_n(\tau)
&\isdef \int_{x=0}^1 e^{2\p inx} \p_\n\matr1x01 g(\tau)
\,dx\\
&= \int_0^1 e^{2\p inx} \left(
\frac{\tau-i}{\tau-x-i}\right)^{(1+\n)/2} \left(
\frac{\tau+i}{\tau-x+i} \right)^{(1+\n)/2}
g(\tau-x)\,dx.
\end{align*}
This defines $g_n$ on a set~$U$ contained in the domain
of~$g$ such that $T\cup U$ is a neighborhood of~$T$.
Let $\Fkr_n\al$ be the hyperfunction represented
by~$g_n$; this does not depend on the choice of the
representative~$g$. One can check that $\p_\n\matr1x01
\Fkr_n\al = e^{-2\p inx} \Fkr_n\al$.

In Proposition~\ref{alphtof} we have discussed the
function~$f_\al$; see also~\vgl{falphdef}. On a
neighborhood of~$T\setminus\vz\infty$ we can replace
$g$ by $\tau\mapsto \left(1+\tau^2\right)^{(1+\n)/2}
f_\al(\tau)$. Then we obtain the following
representative~$ g^\ast_n$ of the restriction
of~$\Fkr_n\al$ to $T\setminus\vz\infty$:
\[ \widearray{r|c|c}
& 0<\im\tau<1& -1<\im\tau<0 \\ \hline
n>0 & g^\ast_n(\tau) = A_n(\al) p_n(\tau) &
g^\ast_n(\tau) = 0 \\
n=0 & g^\ast_0(\tau) = \frac12 A_0(\al) p_0(\tau) &
g^\ast_0(\tau) = -\frac12 A_0(\al) p_0(\tau) \\
n<0 & g^\ast_n(\tau) = 0
& g^\ast_n(\tau) = -A_n(\al) p_n(\tau) \\
\end{array}
\]
Here $p_n(\tau) = \left( 1+\tau^2\right)^{(1+\n)/2}
e^{2\p in\tau}$. For $n\neq0$ we conclude that
$\Fkr_n\al = A_n(\al) \k_n(\n)$.

By interchanging the order of integration we obtain
$\left\langle
\p_{-\n}(p(z))\ph_{2r},\Fkr_n\al\right\rangle =
\int_0^1 e^{2\p inx'} \allowbreak \left\langle
\p_{-\n}(p(z))\ph_{2r},\p_\n\matr1{x'}01\al \right
\rangle \,dx' = F_n u_{2r}(p(z)) $.

\ntst{Notation.}\label{falph0c}We write
$f_\al=f_\al^0+f_\al^c$, with $f_\al^0(\tau) \isdef
\frac12 A_0(\al)\sign(\im\tau)$.

\begin{stel}{Lemma}\label{cuspdecomp}Let $\al \in
A^\n_{-\om}(\G)$, and suppose that it has polynomial
growth at~$\infty$. There is a decomposition
$\al=\Fkr_0\al +\al^c$, with $\al^c\in M^\n_{-\om}$,
such that the function $\tau\mapsto
\left(1+\tau^2\right)^{(1+\n)/2} f_\al^c(\tau)$
represents the restriction of $\al^c$ to
 $T\setminus\vz\infty$, and such that the Fourier term
 of order~$0$ is given by
\[ \Fkr_0\al =
\begin{cases}
B_0(\al) \m + \p^{-1/2}\Gf(-\n/2) A_0(\al) \k_0(\n)
\quad&\text{if $\n\not\in \ZZ_{\geq0}$,}\\
B_0(\al) \m - A_0(\al) \ld(0) \quad&\text{if $\n=0$,}\\
B_0(\al) \m + C_0(\al) \m^\ast_\n + \p^{-1/2}\Gf(-\n/2)
A_0(\al) \k_0(\n) \quad
&\text{if $\n\geq1$ is odd,}\\
B_0(\al) \m + C_0(\al) \k_0(\n) \quad&\text{if
$\n\geq2$ is even.}
\end{cases}
\]
\end{stel}
\bw{Remark}See \ref{intNorbit} and~\ref{kap0cont} for
the definition of the hyperfunction~$\k_0(\n)$,
\ref{muinf} for the definition of~$\m$
and~$\m^\ast_\n$, and~\ref{lambdadef} for~$\ld(0)$.
\bw{Proof}See \ref{alphcdef}--\ref{Ft0}.

\ntst{Definition of~$\al^c$.}\label{alphcdef}The
function $f_\al^c$ decreases quickly at~$\pm i\infty$.
This implies that for each $\ph\in \Akr_T(T)$ the
integrals in
\[ \frac1\p \sum_\pm \int_{I_\pm} \ph(\tau) \left( \pm
f_\al^c(\tau) \right) \left(1+\tau^2\right)^{(\n-1)/2}
\, d\tau\]
converge absolutely, if we take contours $I_+$
and~$I_-$ as indicated in Figure~\ref{Icont} on
page~\pageref{Icont} inside the domain of~$\ph$. In
the same way as in~\ref{intNorbit} we show that this
linear form is given by an hyperfunction. We define
$\al^c$ to be this hyperfunction. The $1$-periodicity
of~$f_\al^c$ gives $\p_\n\matr1101 \al^c=\al^c$.

\ntst{Representative of~$\al^c$.}\label{alphcrepr}Like
we did in~\ref{kapnrepr}, we obtain a representative
of~$\al^c$ by defining
\begin{align}
g_\al^c(\tau_0) &\isdef \frac1{2\p i} \sum_\pm
\int_{I_\pm} \frac{1+\tau\tau_0}{\tau-\tau_0}
\left(\pm f_\al^c(\tau)\right)
\left(1+\tau^2\right)^{(\n-1)/2}\,d\tau
\displaybreak[0] \label{alphcreprint}\\
&= \left(1+\tau_0^2\right)^{(1+\n)/2}
f_\al^c(\tau_0)\nonumber\\
&\quad\hbox{} + \frac1{2\p i} \sum_\pm \int_{I_\pm}
\frac1{\tau-\tau_0} \left( \pm f_\al^c(\tau) \right)
\left(1+\tau^2\right)^{(\n+1)/2}\,d\tau \nonumber \\
&\quad\hbox{} - \frac1{2\p i} \sum_\pm \int_{I_\pm}
\left( \pm f_\al^c(\tau) \right) \tau
\left(1+\tau^2\right)^{(\n-1)/2}\,d\tau.
\label{alphcreprext}
\end{align}
In~\vgl{alphcreprint} the point~$\tau_0\in\hv^\e$, with
$\e=+$ or~$-$, is supposed to be inside the
contour~$I_\e$. By taking the contours wide enough, we
conclude that $g_\al^c(\tau_0) = \oh(1)$ as
$|\im\tau_0|\rightarrow\infty$, uniformly for
$\re\tau_0$ in compact sets.

In~\vgl{alphcreprext} we suppose that $\tau_0\in\hv^\e$
is between the real axis and the contour~$I_\e$. The
integrals in~\vgl{alphcreprext} define functions that
are holomorphic on a neighborhood of~$\RR$. The latter
of the integrals does not depend on~$\tau_0$. The
function $\tau \mapsto
\left(1+\tau^2\right)^{(1+\n)/2} f^c_\al(\tau)$
represents the restriction of~$\al^c$
to~$T\setminus\vz\infty$.

\ntst{Fourier expansion.}Consider
$\ph=\p_{-\n}(p(z))\ph_{2r}$, and insert the series
expansion of~$f_\al^c$. We interchange the order of
summation and integration, and obtain
\begin{equation}
\label{alphcseries}
\left\langle \p_{-\n}(p(z))\ph_{2r},\al^c\right\rangle
= \sum_{n\neq0} e^{2\p inx} \frac{ (-1)^r
(\p|n|)^{-(1+\n)/2} A_n(\al) } {\Gf(\frac{1-\n}2+
r\sign n)} W_{r\sign n,\, \n/2} (4\p |n|y).
\end{equation}

{}From~\ref{Ftal} and the fact that hyperfunctions are
determined by their values on the~$\ph_{2r}$, it
follows that $\al-\al^c=\Fkr_0\al$.

\ntst{The Fourier term of order zero.}\label{Ft0}In
Table~\ref{Fouriersystems} we see that
$\left(F_0u_{2r}\right)_r$ is a linear combination of
systems of Fourier terms. These are represented by
$N$-invariant hyperfunctions, as indicated in
Table~\ref{Neqvarhf}. We have seen that the
restriction of $\Fkr_0\al$ to $T\setminus\vz\infty$ is
represented by $\tau\mapsto \frac12 A_0(\al)
\allowbreak \sign\left(\im\tau\right) \allowbreak
\left(1+\tau^2\right)^{(1+\n)/2}$. In
Table~\ref{supprestr} we check what the multiples of
$\k_0(\n)$, respectively $\ld(0)$, have to be; we see
also that $A_0(\al)$ has to vanish if $\n\geq2$ is
even. Thus we get the expression for~$\Fkr_0\al$ in
the lemma; this serves as the definition of $B_0(\al)$
and~$C_0(\al)$.

\ntst{Fourier coefficients.}\label{fcffs}We have seen
in the course of the proof of Lemma~\ref{cuspdecomp}
that for $\pm n>0$
\[ F_nu_{2r}(p(z)) = \frac{(-1)^r (\p|n|)^{-(1+\n)/2}
}{\Gf(\frac{1-\n}2\pm r) } A_n(\al) e^{2\p inx} W_{\pm
r,\n/2}(4\p|n|y),\]
and that $F_0u_{2r}(p(z)) $ is equal to
\begin{multline*}
B_0(\al) y^{(1-\n)/2} \\
\hbox{} +
\begin{cases}
A_0(\al) \frac{
\Gf(-\frac\n2)\left(\frac{1+\n}2\right)_{|r|} }
{\sqrt\p \Gf(\frac{1-\n}2+|r|)} y^{(1+\n)/2} &\text{
if $\n\in\CC\setminus\ZZ_{\geq0}$,}\\
\frac{-2}\p A_0(\al) y^{1/2}\left( \log y + l_r\right)
&\text{ if $\n=0$,}\\
\left( A_0(\al) -i \sign(r) C_0(\al) \right)\frac{
\Gf(-\frac\n2)\left(\frac{1+\n}2\right)_{|r|} }
{\sqrt\p \Gf(\frac{1-\n}2+|r|)} y^{(1+\n)/2}
 &\text{ if $\n\geq1$ is odd,}\\
C_0(\al) \frac{ \left(\frac{\n+1}2\right)_{|r|} }{
\Gf(\frac{1-\n}2+|r|)} y^{(1+\n)/2} &\text{ if
$\n\geq2$ is even.}
\end{cases}
\end{multline*}

\ntst{Isomorphism.}\label{isoAB}Let
$\n\in\CC\setminus\left(1+2\ZZ\right)$. The
isomorphism $\iota(\n):M^\n_{-\om} \rightarrow
M^{-\n}_{-\om}$ induces an isomorphism
$\iota(\n):A^\n_{-\om}(\G) \rightarrow
A^{-\n}_{-\om}(\G)$, see~\vgl{iotadef}. It preserves
polynomial growth at~$\infty$. A consideration of
Fourier terms leads to
\begin{alignat*}2 A_n(\iota(\n)\al) &= (\p|n|)^{-\n}
\Gf\left(\frac{1+\n}2\right) \Gf\left(
\frac{1-\n}2\right)^{-1} A_n(\al)& \quad&\text{if
$n\neq0$,}\\
B_0(\iota(\n)\al) &= \frac1{\sqrt\p} \Gf\left( -
\frac\n2\right) \Gf\left( \frac{1-\n}2\right)^{-1}
A_0(\al) &\quad&\text{if $\n\not\in\ZZ_{\geq0}$,}\\
&= \Gf\left(\frac{1+\n}2\right)^{-1} C_0(\al )
&\quad&\text{for $\n\in 2\NN$,}\\
A_0(\iota(\n)\al) &= \sqrt\p
\Gf\left(\frac{1+\n}2\right) \Gf\left(
\frac\n2\right)^{-1} B_0(\al) &\quad&\text{if
$\n\not\in\ZZ_{\leq0}$,}\\
C_0(\iota(\n)\al) &= \Gf\left(\frac{1+\n}2\right)
B_0(\al)
 &\quad&\text{for $\n\in -2\NN$.}
\end{alignat*}
The isomorphism $\iota(0)$ is the identity.

\ntst{Maass forms.}For the hyperfunction~$\al\in
A^{2s-1}_{-\om}(\Gmod)$ associated to the cuspidal
Maass form in~\vgl{Maassform} we have
$A_0(\al)=B_0(\al)=0$ and $A_n(\al) =
\left(\p|n|\right)^s \Gf\left(1-s\right) a_n$ for
$n\neq0$. We have chosen $\n=2s-1$; the choice
$\n=1-2s$ would be as good; it gives
$\iota\left(2s-1\right)\al \in
A^{1-2s}_{-\om}(\Gmod)$.

\ntst{Eisenstein series.}\label{EisAB}In~\ref{Eishf}
we have given $\e^\ast_s \in A^{1-2s}_{-\om}(\Gmod)$
for $\re s>1$. Here the choice $\n=1-2s$ seems the
natural one if one tries to get a hyperfunction from
the series in~\vgl{Eisseries}. {}From the Fourier
expansion in~\vgl{EisFourier} we obtain, for $s
\not\in 1+\frac12\ZZ_{\leq0}$:
\begin{alignat*}2 A_n(\e^\ast_s) &= 2\p^{1-s}
\Gf(s)\s_{1-2s}(|n|)&\quad&\text{for $n\neq0$}
\displaybreak[0] \\
A_0(\e^\ast_s)&= 2\sqrt\p \frac{ \Gf(s) }{
\Gf(s-\frac12) } \Ld(2s-1)
&& \displaybreak[0] \\
B_0(\e^\ast_s) &= 2\Ld(2s) &&
\end{alignat*}

This suggests to consider the family $\e_s\isdef
\frac12\G\left(s-1\right)^{-1} \e^\ast_s$.
\begin{stel}{Proposition}\label{Eisfam}The family
$s\mapsto \e_s$ is holomorphic on~$\CC$.

$\e_s\in A^{1-2s}_{-\om}(\Gmod)$ for each $s\in\CC$.
\end{stel}
\bw{Remark}We call a family $s\mapsto \e_s$ of
hyperfunctions holomorphic if $s\mapsto \left\langle
\ph,\e_s\right\rangle$ is holomorphic for each $\ph\in
\Akr_T(T)$.
\bw{Proof}For $\tau\in\hv^\pm$ we have
\begin{equation}
\label{fceis}
f_{\e_s}^c(\tau) = \pm \p^{1-s} \left(s-1\right)
\sum_{n=1}^\infty \s_{1-2s}(n) e^{\pm 2 \p in\tau}.
\end{equation}
This converges uniformly for~$\tau$ in compact sets.
Thus we obtain holomorphy of~$s\mapsto \e_s^c$.

The Fourier term of order~$0$ is
\begin{alignat}2 \label{F0eis}
\Fkr_0\e_s &= \frac{\Ld(2s)}{\Gf(s-1)} \m +
(s-1)\Ld(2s-1) \k_0(1-2s)&\quad&\text{for
$s\neq\frac12$,} \\
\label{F0eis0}
&= \left( \frac{\Ld(2s)}{\Gf(s-1)} + \p^{-1/2} (s-1)
\Ld(2s-1) \right) \m \\
\nonumber
&\qquad- \sqrt\p \left(s-\frac12\right) (s-1) \Ld(2s-1)
\ld(1-2s)&\quad&\text{for $s$ near $\frac12$.}
\end{alignat}
This shows the holomorphy of~$s\mapsto \e_s$. The
$\Gmod$-invariance extends by holomorphy.

\ntst{Functional equation.} We have $\Gf(-s) \e_{1-s} =
\Gf\left(1-s\right) \iota\left(2s-1\right) \e_s$ for
$s\not\in\ZZ$.

\nwsect{Geodesic decomposition}\label{secgeodec}
\ntst{Period polynomials.}\label{clcocycle}Let $H$ be a
holomorphic modular cusp form of weight~$2k\geq12$.
The period polynomial~$r_H$ associated to~$H$ is
\[ r_H(X) \isdef \int_0^\infty H(\tau) (X-\tau)^{2k-2}
\, d\tau,\]
see, e.g., \cite{Za}. The cusps $0$ and~$\infty$ can be
replaced by any pair $(\x,\eta)$ of cusps. The path of
integration should approach $\x$ and~$\eta$ along a
geodesic for the non-Euclidean metric on~$\bhv$. In
this way we arrive at a homogeneous $1$-cocycle~$R_H$
with values in the polynomials of degree at
most~$2k-2$:
\[ R_H(\x,\eta;X) \isdef \int_\x^\eta H(\tau)
(X-\tau)^{2k-2}\,d\tau. \]
It satisfies $R_H\left(\g\cdot\x,\g\cdot\eta\right) =
R_H(\x,\eta)|_{2-2k}\g^{-1}$ for $\x,\eta,\th\in\prq$,
$\g\in\Gmod$. (See~\ref{polact} for the action
$F\mapsto F_{2-2k}g$.)

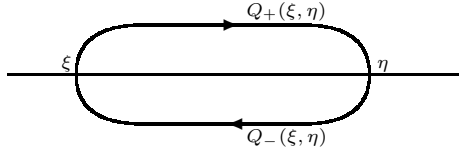
\begin{figure}
\begin{center}
\scriptsize
\setlength\unitlength{.6cm}
\begin{picture}(10,2.5)(-5,-1.2)
\put(-5,0){\line(1,0){10}}
\put(-3.8,.15){$\x$}
\put(3.2,.15){$\eta$}
\put(.3,1.3){$Q_+(\x,\eta)$}
\put(.3,-1.5){$Q_-(\x,\eta)$}
\thicklines
\put(-.15,1.1){\vector(1,0){.3}}
\put(.15,-1.1){\vector(-1,0){.3}}
\qbezier(-3.5,0)(-3.5,1.1)(-2,1.1)
\qbezier(3,0)(3,1.1)(1.5,1.1)
\qbezier(-3.5,0)(-3.5,-1.1)(-2,-1.1)
\qbezier(3,0)(3,-1.1)(1.5,-1.1)
\put(-2,1.1){\line(1,0){3.5}}
\put(-2,-1.1){\line(1,0){3.5}}
\end{picture}
\end{center}
\caption[o]{The contour $Q(\x,\eta)$ used in the
definition of the geodesic decomposition is the union
of $Q_+(\x,\eta)$ and $Q_-(\x,\eta)$. Near $\x$
and~$\eta$ the contours are pieces of conjugate
geodesics in $\bhv$ and~$\ohv$.} \label{Qcont}
\end{figure}

\ntst{Discussion.}We want to generalize this to
hyperfunctions associated to cuspidal Maass forms and,
as far as possible, to Eisenstein series. We try to
integrate a representative~$g$ of the automorphic
hyperfunction~$\al$ along a path as given in
Figure~\ref{Qcont} on page~\pageref{Qcont}. If $g$
stays bounded on geodesics approaching the cusps $\x$
and~$\eta$, then this is no problem. This holds for
$\x,\eta\in \G\cdot\infty$ if $\al=\al^c$. But if we
use a principal value interpretation of the integral
near~$\x$ and~$\eta$, we can extend this approach to
more automorphic hyperfunctions.

We arrive at quantities $\al[\x,\eta] \in
M^{\n}_{-\om}$ that just fail to be cocycles. In the
case that $\al$ corresponds to a holomorphic cusp
form, the map $M^{2k-1}\rightarrow E(2k)$
in~\ref{topol} sends $\al[\x,\eta]$ to a multiple
of~$R_H(\x,\eta)$.

\ntst{Definition.}\label{Gamdecompdef}Let $X\subset T$,
$X\neq\emptyset$, be invariant under~$\G$. We call a
{\sl $\G$-decomposition}~$p$ of $\al\in
A^\n_{-\om}(\G)$ {\sl on~$X$}\/ a map
$\vzm{(\x,\eta)\in X^2}{\x\neq\eta} \rightarrow
M^\n_{-\om}: (\x,\eta) \mapsto \al[\x,\eta]_p$ that
satisfies
\begin{enumerate}
\item[a)] $\supp(\al[\x,\eta]_p) \subset [\x,\eta]$ for
all $\x,\eta\in X$, $\x\neq\eta$.
\item[b)] $\al=\sum_{j=1}^n \al[\x_{j-1},\x_j]_p = \al$
whenever $\x_1,\ldots,\x_n\in X$ satisfy $\x_0<\x_1<
\cdots<\x_n=\x_0$.
\item[c)] $\al[\g\cdot\x,\g\cdot\eta] =
\p_\n(\g)\al[\x,\eta]$ for all $\g\in\G$ and
$\x,\eta\in X$, $\x\neq\eta$.
\end{enumerate}
In condition~a) the closed interval $[\x,\eta]\subset
T$ is understood to refer to the cyclic order on the
circle~$T$. In condition~b) the~$\x_j$ are supposed to
go around~$T$ only once: the intervals
$[\x_{j-1},\x_j]$ intersect each other only in the end
points.

We define $A^\n_{-\om}(\G,X) = \vzm{\al\in
A^\n_{-\om}(\G)}{\text{$\al$ has a $\G$-decomposition
on~$X$}}$.

\ntst{Partings and $\G$-decompositions.}Let $\x\in X$,
and let $I\neq T$ an open interval containing~$\x$
such that $T\setminus I$ contains a point ~$\eta \in
X$. The image of the decomposition $\al =
\al[\eta,\x]_p + \al[\x,\eta]_p$ in $\Bkr_T(I)$
determines a parting of~$\al$ at~$t$
(see~\ref{decompR}). The map $\RR\rightarrow
T:\th\mapsto\cot\th$ that we used to define~$\Bkr_T$
is decreasing. So for a parting at~$\x$ of a
hyperfunction on~$T$ the support of~$\al_-$ is to the
left of~$\x$, and that of~$\al_+$ to the right.

Conversely, if we have a parting $\al = \al_{\x,-} +
\al_{\x,+}$ in the stalk~$\left( \Bkr_T\right)_\x$ for
each $\x\in X$, we get a decomposition satisfying a)
and~b) in~\ref{Gamdecompdef}: Take open intervals
$I_\x$ and~$I_\eta$ containing~$\x$,
respectively~$\eta$, with empty intersection and
determine $\al[\x,\eta]_p$ by its restrictions:
\begin{alignat*}2 \left. \al[\x,\eta]_p \right|_{I_\x}
&= \al_{\x,+}
&\quad&\text{on $ I_\x $,}\\
\left. \al[\x,\eta]_p \right|_{(\x,\eta)} &= \al
&\quad&\text{on $ (\x,\eta) $,}\\
\left. \al[\x,\eta]_p \right|_{I_\eta} &= \al_{\eta,-}
&\quad&\text{on $ I_\eta $,}\\
\left. \al[\x,\eta]_p \right|_{(\eta,\x)} &= 0
&\quad&\text{on $(\eta,\x) $.}
\end{alignat*}
Condition~c) is equivalent to $\p_\n(\g)\al_{\x,\pm} =
\al_{\g\cdot\x,\pm}$ for all $\x\in X$, $\g\in\G$.

Let us write $X$ as a disjoint union of $\G$-orbits:
$X= \bigsqcup_{\x\in\X} \G\cdot\x$. Finding a
$\G$-decomposition of~$\al$ on~$X$ is equivalent to
finding a parting $\al=\al_{\x,-} + \al_{\x,+}$ at
each $\x\in\X$ such that $\n_\n(\dt) \al_{\x,\pm} =
\al_{\x,\pm}$ for each $\dt\in \G_\x \isdef
\vzm{\g\in\G}{\g\cdot\x=\x}$. Thus we have obtained:

\begin{stel}{Proposition}Let $\al \in A^\n_{-\om}(\G)$.
\begin{enumerate}
\item[i)] Let $\x\in T$. Then $\al \in
A^\n_{-\om}\left(\G,\G\cdot\x\right)$ in the following
cases:
\begin{enumerate}
\item[a)] $\G_\x=\vz1$.
\item[b)] There is a parting $\al=\al_-+\al_+$ of~$\al$
at~$\x$ that satisfies $\p_\n(\dt)\al_\pm=\al_\pm$ in
the stalk~$\left(\Bkr_T\right)_\x$ for all
$\dt\in\G_\x$.
\end{enumerate}
\item[ii)] Let $\vz{X_j}_{j\in J}$ be a collection of
non-empty, $\G$-invariant subsets of~$T$, and put
$X\isdef \bigcup_{j\in J} X_j$. Then $\al \in
A^\n_{-\om}\left(\G,X\right)$ if and only if $\al \in
A^\n_{-\om}\left(\G, X_j\right)$ for all $j\in J$.
\end{enumerate}
\end{stel}
\bw{Remark}It may very well be true that
$A^\n_{-\om}(\G,X) = A^\n_{-\om}(\G)$ for all~$X$. In
this paper we direct our attention to $X=\G\cdot
\infty$.

\ntst{Definition.}\label{geodappr}Let $\al$ be a
hyperfunction (not necessarily automorphic),
represented by $g\in\Okr\left( U\setminus T\right)$
for some neighborhood~$U$ of~$T$ in~$\prc$. We define
$\al$ to have {\sl geodesic approach at~$\infty$} if
for each holomorphic function~$\ph$ on a neighborhood
of~$\infty$ the following conditions are satisfied:
\begin{enumerate}
\item[a)] For each $\x\in\RR$ and for each sufficiently
large~$y>0$ the function
\[ t\mapsto \frac i\p \left(
\frac{\ph(\x+it)g(\x+it)}{1+(\x+it)^2} +
\frac{\ph(\x-it)g(\x-it)}{1+(\x-it)^2} \right) \]
is integrable on~$[y,\infty)$ (with respect to the
measure~$dt$).
\item[b)] For all $\x,\x_1\in\RR$:
\[ \lim_{Y\rightarrow\infty} \int_{x=\x}^{\x_1} \frac
1\p \left( \frac{ \ph(x+iY) g(x+iY)} {1+(x+iY)^2} -
\frac{ \ph(x-iY) g(x-iY)} {1+(x-iY)^2} \right)
\,dx=0.\]
\end{enumerate}

Let $m\in G$. We define $\al$ to have geodesic approach
at $m \cdot \infty$ if $\p_\n(m )\al$ has geodesic
approach at~$\infty$. This does not depend on the
choice of~$\n\in\CC$.

\ntst{Discussion.}The integral of the function in
condition~a) is a principal value variant of the
integral $ \int_{L_\x} \ph(\tau)g(\tau)
\frac{d\tau}{\p(1+\tau^2)}$, where $L_\x$ is the path
from $\x+iy$ vertically upward to~$\infty$ in~$\bhv$,
and then from~$\infty$ vertically upward in~$\ohv$ to
$\x-iy$. If $g$ is bounded on a neighborhood
of~$\infty$, then this integral exists. The
formulation in condition~a) allows some cancellation
between both parts of the integral. Similarly, the
integral in condition~b) is a principal value form of
$ \int_{M(\x,\x_1;Y)} \ph(\tau)g(\tau)
\frac{d\tau}{\p(1+\tau^2)}$ where $M(\x,\x_1;Y)$
consists of a path from $\x+iY$ to $\x_1+iY$ and a
path from $\x_1-iY$ to $\x-iY$. These two conditions
allow us to define integrals of
$\ph(\tau)g(\tau)\frac{d\tau}{\p(1+\tau^2)}$ from a
point $z\in\bhv$ via~$\infty$ to the points~$\bar
z\in\ohv$, where $\infty$ is crossed along two
conjugate geodesics in $\bhv$ and~$\ohv$. The
particular choice of the geodesic does not matter. The
choice of the representative~$g$ does not influence
the conditions.

Let $p\in G$ be such that $p\cdot\infty=\infty$. Then
$\p_\n (p)\al$ has a representative $g_1:\tau\mapsto
J(\tau)^{(1+\n)/2}g\left( p^{-1} \cdot\tau\right)$,
with $J$ holomorphic on a neighborhood of~$\infty$. As
$p$ has the form $\matr \ast\ast0\ast$, conditions a)
and~b) for $g$ and~$g_1$ are equivalent. This implies
that the definition of geodesic approach at $
m\cdot\infty $ does not depend on~$m$, only on $
m\cdot\infty $.

\ntst{Principal value integrals.}\label{pv}Geodesic
approach at $\x\in T$ allows us to integrate over
paths passing the point~$\x$ along any pair of
conjugate geodesics we like.

If $L$ is a path that crosses~$T$ at a finite number of
points along pairs of conjugate geodesics, we denote
by $\pvint_L \ph(\tau) g(\tau)
\frac{d\tau}{\p(1+\tau^2)}$ the integral over~$L$ in
which the contributions along each of the pairs of
geodesics has the interpretation given above.

If we deform~$L$ in such a way that the end points and
the points at which $T$ is crossed are kept fixed (the
corresponding pairs of conjugate geodesics may
change), then the integral does not change.

\begin{stel}{Lemma}\label{Fplusminlemma}Let the
function~$g$ represent $\al\in M^\n_{-\om}$ on a
neighborhood of~$\infty$. Define for $\x\in\RR$ and
$t$ large:
\begin{align*}
F_+(\x,t) &\isdef \frac12\left( \frac{ g(\x+it)
}{1+(\x+it)^2} + \frac{g(\x-it)} {1+(\x-it)^2}
\right),\\
F_-(\x,t) &\isdef \frac1{2t} \left( \frac{ g(\x+it)
}{1+(\x+it)^2} - \frac{g(\x-it)} {1+(\x-it)^2}
\right).
\end{align*}

\begin{enumerate}
\item[i)] Condition~a) in the definition of geodesic
approach at~$\infty$ in~\ref{geodappr} is equivalent
to the integrability (for the measure~$dt$) of all
$F_\pm(\x,\cdot)$ on each interval $[y,\infty)$ with
$y$ large.

\item[ii)] If $F_+(x,t)=o(t)$ and $F_-(x,t)=o(1/t)$ as
$t\rightarrow\infty$ uniformly for $x$ between $\x$
and~$\x_1$, then condition~b) in the definition of
geodesic approach at~$\infty$ is satisfied.

\item[iii)] Let $L_\x$ be the path along from $\x+iy$
via~$\infty$ to $\x-iy$ indicated above, and let $\ph$
be holomorphic on a neighborhood of~$\infty$
containing~$L_\x$. If $\al$ has geodesic approach
at~$\infty$, then $\pvint_{L_\x} \ph(\tau) g(\tau)
\frac{d\tau}{\p(1+\tau^2)}$ is equal to
\[ \frac i\p \int_y^\infty F_+(\x,t) \sum_\pm
\ph\left(\x\pm it\right)\,dt + \frac i\p \int_y^\infty
F_-(\x,t) \sum_\pm \pm t\, \ph\left(\x\pm
it\right)\,dt.\]
\end{enumerate}
\end{stel}
\bw{Proof}Take $\ph(\tau) =1$ and $\ph(\tau) =
\frac1{\tau-\x}$ in condition~a) to see that
$F_{\x,\pm}$ is integrable.

We note that $\sum_\pm \ph\left(\x\pm it\right)$ and
$\sum_\pm \pm t\, \ph\left(\x\pm it\right)$ are
bounded as $t\rightarrow\infty$. (Use the holomorphy
of~$\ph$ at~$\infty$.) A computation shows that the
integral of the function in condition~a) equals the
sums of the integrals in part~iii). This gives
part~iii) and the converse implication in part~i).

Part~ii) follows from the facts that the integral in
condition~b) is equal to
\begin{equation}
\label{bintl}
\frac1\p \int_\x^{\x_1} \left( F_+(x,Y) \sum_\pm \pm
\ph\left(\x\pm iY\right) + F_-(x,Y) \sum_\pm Y
\ph\left(\x\pm iY\right) \right)\,dx ,\end{equation}
and that $\sum_\pm \pm \ph\left(\x\pm iY\right) =
\oh(1/Y)$ as $Y\rightarrow\infty$.

\begin{stel}{Lemma}\label{propgeodappr}The following
hyperfunctions have geodesic approach at~$\infty$:
\begin{enumerate}
\item[i)] $\al^c$ for each $\al\in A^\n_{-\om}(\G)$
with polynomial growth at~$\infty$.
\item[ii)] $\m$ and $\ld(0)$.
\item[iii)] $\k_0(\n)$ for $\re\n<1$.
\end{enumerate}
\end{stel}
\bw{Proof}If a representative~$g$ is bounded on
vertical lines uniformly for $\re \tau$ in compact
sets, then the functions~$F_\pm$ in
Lemma~\ref{Fplusminlemma} satisfy $F_\pm(\x,t) =
\oh\left(t^{-2}\right)$. This suffices for~$\al^c$,
see~\ref{alphcrepr}.

The representative $g(\tau) = \frac{-i}2\tau$ of~$\m$
satisfies $F_\pm(\x,t) = \oh(t^{-2})$.

For part~iii) we use the representative~$q$
in~\ref{kap0cont} for $\n\not\in -2\NN$. This gives
$F_\pm(\x,t) = \oh\left(t^{\re \n-2}\right)$. For
$\n\in-2\NN$ we proceed as for~$\al^c$.

Finally, we check that $F_\pm(\x,t) =
\oh\left(t^{-2}\log t \right)$ for~$\ld(0)$.

\begin{stel}{Lemma}\label{corgd}Let $\al \in
A^\n_{-\om}(\G)$ have polynomial growth at~$\infty$.
Then $\al$ has geodesic approach at all points of the
orbit~$\G\cdot\infty$ if one of the following
conditions is satisfied:
\begin{enumerate}
\item[a)] $\re \n<1$.
\item[b)] $A_0(\al)=0$, and if $\n\in\NN$, then
$C_0(\al)=0$.
\end{enumerate}
\end{stel}
\bw{Proof}The geodesic approach at~$\infty$ follows
directly from Lemma~\ref{propgeodappr}. Use the
$\G$-invariance for the other points
of~$\G\cdot\infty$.%

\ntst{Geodesic decomposition.}Let $X\subset$ be the set
of points at which a given hyperfunction~$\bt$ has
geodesic approach.

Let $\x,\eta\in X$, $\x\neq\eta$. We define for $\ph\in
\Akr_T(T)$:
\[ \left\langle \ph,\bt[\x,\eta]_g \right\rangle \isdef
\pvint_{Q(\x,\eta)} \ph(\tau) g(\tau)
\frac{d\tau}{\p(1+\tau^2)}, \]
with the contour~$Q(\x,\eta)$ given in
Figure~\ref{Qcont}. It is understood that the region
between $Q(\x,\eta)$ and the interval $[\x,\eta]$ is
contained in the domain of~$\ph$. To see that this
defines $\bt[\x,\eta]_g$ as a hyperfunction, we
estimate $\left\langle \ph_{2r},\bt[\x,\eta]_g
\right\rangle$ in terms of the supremum norm of
$\ph_{2r}$ and its first derivative (with respect to a
local coordinate). (Use part~iii) of
Lemma~\ref{Fplusminlemma}.)

Consider $\left\langle h_{\tau_0},\bt[\x,\eta]_g
\right\rangle$ (see~\ref{repres}), with $\tau_0$
outside $Q(\x,\eta)$, to see that the hyperfunction
$\bt[\x,\eta]_g$ has support contained in $[\x,\eta]$.
If $\th\in X \cap (\x,\eta)$, then $\bt[\x,\eta]_g =
\bt[\x,\th]_g + \bt[\th,\eta]_g$.

If $T=\bigcup_{j=1}^m [\x_j,\x_{j+1}]$ is a partition
of~$T$ with $\x_{m+1}=\x_1 < \cdots < \x_m\in X$ such
that the intervals $[\x_j,\x_{j+1}]$ intersect each
other only in their end points, then $\bt =
\sum_{j=1}^m \bt[\x_j,\x_{j+1}]_g$. This we call the
{\sl geodesic decomposition on~$X$}. In the sequel we
write $\bt[\x,\eta]$ instead of $\bt[\x,\eta]_g$.

\ntst{$\G$-behavior.}Let $m\in G$. If $\x,\eta,
m\cdot\x, m\cdot\eta\in X$, then a computation shows
that $\p_\n(m) \left( \bt[\x,\eta]\right) = \left(
\p_\n(m)\bt\right)\left[m\cdot\x,m\cdot\eta\right]$.
So, if $\al$ is an automorphic hyperfunction, the
set~$X$ of points at which it has geodesic approach is
$\G$-invariant, and the geodesic decomposition
of~$\al$ on~$X$ is a $\G$-decomposition.

\ntst{Reflection.}If $\bt$ has geodesic approach at
$\x,\eta,-\x,-\eta$, then $\spgl\left(
\bt[\x,\eta]\right) = \left( \spgl\al\right)
[-\eta,-\x]$.

\begin{stel}{Theorem}\label{thmgdc}Let $\al\in
A^\n_{-\om}(\G)$ have polynomial growth at~$\infty$.
Then $\al$ has geodesic decomposition on the
orbit~$\G\cdot \infty$ if one of the following
conditions is satisfied:
\begin{enumerate}
\item[a)] $\re \n<1$.
\item[b)] $A_0(\al)=0$, and if $\n\in\NN$, then
$C_0(\al)=0$.
\end{enumerate}
\end{stel}
\bw{Remarks}One finds the Fourier coefficients
$A_0(\al)$ and $C_0(\al)$ in Lemma~\ref{cuspdecomp}.

If $\al$ corresponds to a system cusp forms, then
condition~b) holds.

Of course, if $\G$ has more than one cuspidal orbit, we
may move $\G$-inequivalent cusps to~$\infty$ by
conjugation. In particular if at each cusp one of the
conditions holds, then we have geodesic decomposition
on the set of all cusps. Note that the theorem allows
$\al$ to have terrible growth at cusps that are not
$\G$-equivalent to~$\infty$.

\bw{Proof}Directly from Lemma~\ref{corgd}.

\ntst{Holomorphic cusp forms.}\label{perfholcf}Let $\al
\in H^{2k-1}_{-\om}(\G)$ correspond to a holomorphic
cusp form $H$ for~$\G$ of weight~$2k\geq2$
(see~\ref{hypholo}). Then $\al$ has geodesic
decomposition on the set~$X$ of all $\G$-cusps. We
consider the image under the map $\bt\mapsto
\bt^{\langle 2k \rangle}$ discussed in~\ref{topol}.
For $\x,\eta\in X$, $\x\neq\eta$, we find the period
integral discussed in~\ref{clcocycle}:
\begin{align*}
\left( \al[\x,\eta]\right)^{\langle 2k \rangle} &=
\int_{Q_+(\x,\eta)} h_X(\tau) (-1)^k 4^{-k}
\left(1+\tau^2\right)^{k} H(\tau)
\frac{d\tau}{\p(1+\tau^2)} \displaybreak[0] \\
&= \frac1{4\p} \int_\x^\eta H(\tau)
\left(\tau-X\right)^{2k-2}\,d\tau = \frac1{4\p}
R_H(\x,\eta).
\end{align*}

\ntst{Eisenstein family.}\label{Eisdecomp}The
Eisenstein family~$\e_s$, introduced in~\ref{EisAB},
has geodesic decomposition on~$\prq$ for $\re s>0$.
{}From the explicit expression in~\ref{Eishf} for
$\left\langle\ph, \e^\ast_s\right\rangle$ converging
for $\re s>1$, one expects that the geodesic
decomposition has the following form
\begin{equation}
\label{ed}
 \left\langle \ph, \e_s[\x,\eta] \right\rangle =
 \frac12\p^{-1} (s-1) \sideset{}' \sum_{p,q} f_{p,q}
 \left(p^2+q^2\right)^{-s} \ph\left(-\frac pq \right)
\end{equation}
for $\re s>1$ and $\x,\eta\in\prq$, with $f_{p,q}\isdef
1 $ if $-\frac pq \in (\x,\eta)$, $f_{p,q} \isdef
\frac12$ if $-\frac pq = \x$~or~$\eta$, and $f_{p,q}
\isdef 0$ otherwise.

This turns out to be true. In the proof one uses the
representative $g_s(\tau) = \frac{-i}4 \p^{-s}
\left(s-1\right) \allowbreak \sideset{}'\sum_{p,q}
\left( p^2+q^2\right)^{-s} \allowbreak
\frac{p\tau-q}{q\tau+p}$ in the definition of
$\left\langle \ph, \e_s[\x,\eta] \right\rangle$. The
main point is to interchange the order of integration
and summation. This is no problem on most of
$Q(\x,\eta)$. But near~$\x$, and near~$\eta$, we have
to treat together terms for which the points~$-\frac
pq$ are symmetrical with respect to~$\x$,
respectively~$\eta$.

One can also show that $s\mapsto \left\langle \ph,
\e_s[\x,\eta] \right\rangle $ is holomorphic on $\re
s>0$.

\nwsect{Infinite sums}\label{secadd}The geodesic
decomposition allows us to write some automorphic
hyperfunctions as a finite sum. Under additional
conditions an infinite decomposition is possible.

\ntst{Weak convergence.}Limits and infinite series of
hyperfunctions we consider in the weak sense:
$\lim_{n\rightarrow\infty} \bt_n=\bt$ means
$\lim_{n\rightarrow\infty} \left\langle \ph,
\bt_n\right\rangle = \left\langle\ph, \bt
\right\rangle$ for each $\ph$ that is holomorphic on a
neighborhood in~$\prc$ of $\supp(\bt) \cup
\bigcup_{n\geq N}\supp(\bt_n)$ for some~$N$, and
similarly for convergence of series.

\begin{stel}{Theorem}\label{add}Let $\x\in T$,
$\x\neq\infty$. Suppose that $\al \in A^\n_{-\om}(\G)$
has polynomial growth at~$\infty$, has geodesic
approach at~$\x$, and satisfies the following
conditions:
\begin{enumerate}
\item[a)] $\re\n <1$,
\item[b)] if $A_0(\al)\neq 0$, then $\re\n<0$.
\end{enumerate}
Then
\begin{align*}\al[\x,\infty]
&=\frac12 B_0(\al) \m + \sum_{n=1}^\infty
\al[\x+n-1,\x+n],\\
\al[\infty,\x] &= \frac12 B_0(\al) \m +
\sum_{n=1}^\infty \al[\x-n,\x-n+1].
\end{align*}
\end{stel}
\bw{Remarks}The conditions imply that $\al$ has
geodesic approach at~$\infty$. The $\G$-orbit of~$\x$
contains the set $\x+\ZZ$. So $\al$ has also geodesic
approach at all points $\x+n$, $n\in\ZZ$.

The theorem applies to hyperfunctions associated to
cuspidal Maass forms; then $B_0(\al)=0$. The statement
is in general false for hyperfunctions associated to
holomorphic cusp forms. The hyperfunction~$\e_s$
satisfies the conditions for $\re s>\frac12$. The
presence of the term $\frac12 B_0(\al)\m$ is very
clear if one considers \vgl{ed} on page~\pageref{ed}
for $\re s>1$.

\bw{Proof}See \ref{strprfadd}--\ref{endaddprf}.

\ntst{Reformulation.}\label{strprfadd}
By conjugation of~$\G$ we can arrange that $\x=0$.
Application of the result for $[0,\infty]$
to~$\spgl\al \in A^\n_{-\om}(j(\G))$ gives the result
for~$[\infty,0]$. The conditions stay valid under both
transformations.

Let $\ph$ be holomorphic in a neighborhood of~$\infty$.
The weak interpretation of the series means that we
have to prove $\lim_{n\rightarrow\infty}
\left\langle\ph, \al[n,\infty] \right\rangle= \frac12
B_0(\al)$. We write $\ph(\tau)=\ps(1/\tau)$, with
$\ps$ holomorphic on a compact neighborhood of~$0$. By
 $\|\cdot\|$ we denote the supremum norm on this
 neighborhood. Let $C_n$ be the vertical line
 $\re\tau=n$, and let $g$ denote a representative
 of~$\al$. If $n$ is large enough (depending on~$\ph$
and~$g$) we have $\left\langle \ph, \al[n,\infty]
\right\rangle = I_n(g)$, where $ I_n(h) \isdef
\frac1\p \pvint_{C_n} \ps\left(\frac1\tau\right)
h(\tau) \frac{d\tau}{1+\tau^2}$. Taking the principal
value interpretation into account, we have
\[ I_n(h) = \frac i\p \int_0^\infty \sum_\pm \ps\left(
\frac1{n\pm it}\right) h(n\pm it) \frac{dt}{1+(n\pm
it)^2}. \]

The conditions allow us to write the representative as
$g(\tau) = g_\al^c(\tau) + \p^{-1/2}\G(-\frac\n2)
A_0(\n) g_\k + B_0(\al) g_\m$, with $g_\k$ and $g_\m$
 representatives of~$\k_0(\n)$ and~$\m$. We want to
 show that $I_n(h)=o(1)$ as $n\rightarrow\infty$ for
$h=g_\al^c$, and for $h= g_\k$ if $\re \n<0$, and that
$\lim_{n\rightarrow\infty} I_n(g_\m)=\frac12 \ps(0)$.

\ntst{Contribution of $\al^c$.}Let $\re\n<1$. We use
the expression for~$g_\al^c$ in~\vgl{alphcreprext} on
page~\pageref{alphcreprext}. So $I_n(g_\al^c)$ splits
up into three terms.

The third term in~\vgl{alphcreprext} contributes a
constant to~$g_\al^c$. The corresponding integral is
estimated by
\[ \int_0^\infty \|\ps\| \frac {dt}{|1+(n+it)^2|} \ll
\frac 1n \int_0^\infty \frac{dt}{|(1+it)^2+1/n^2|} =
o(1).\]

The middle term in~\vgl{alphcreprext} is
$\oh(\frac1n)$. The corresponding integral is
$\oh(n^{-2})$.

For the first term we use the periodicity of~$f_\al^c$
and obtain
\[ I_n = \frac i\p \int_0^\infty \sum_\pm
\ps\left(\frac1{n\pm it}\right) f_\al^c(0\pm it)
\left(1+(n\pm it)^2\right)^{(\n-1)/2} \,dt.\]
The integral over $[1,\infty)$ we estimate by
$\int_1^\infty \|\ps\| e^{-\e t} \,
\left|n+it\right|^{\re\n-1}dt$. This is $ \oh\left(
n^{\re\n-1}\right) =o(1) $, for some $\e>0$.

The function $\tau\mapsto f_\al^c(\tau)
\left(1+\tau^2\right)^{(1+\n)/2}$ represents the
restriction of~$\al$ to $T\setminus\vz\infty$. As
$\al$ has geodesic approach at~$n$, the integral
$\int_0^1 \sum_\pm \ph_1(n\pm it) \allowbreak
f_\al^c(\pm it) \allowbreak \left(n+(1\pm
it)^2\right)^{(\n-1)/2} dt$ converges for all~$\ph_1$
holomorphic on a neighborhood of $n+i[-1,1]$. Put
$G_+(t) \isdef f_\al^c(\pm it)$, and $G_-(t) \isdef it
\sum_\pm \pm f_\al^c(\pm it)$. The choice
$\ph_1(\tau)=\left(1+\tau^2\right)^{(1-\n)/2}$ shows
that $G_+$ is integrable on~$[0,1]$. The integrability
of~$G_-$ follows from $\ph_1(\tau) =
\left(\tau-n\right)\allowbreak
\left(1+\tau^2\right)^{(1-\n)/2}$. We obtain
\begin{align*}
\int_0^1 \sum_\pm
&\ps\left( \frac1{n\pm it}\right) f_\al^c(\pm it)
\left( 1+(n\pm it)^2 \right)^{(\n-1)/2}\, dt\\
&= \int_0^1 \frac12 G_+(t) \sum_\pm \ps\left(
\frac1{n\pm it}\right) \left( 1+(n\pm it)^2
\right)^{(\n-1)/2}\, dt\\
&\quad\hbox{} + \int_0^1 \frac1{2i} G_-(t) \frac 1t
\sum_\pm \pm \ps\left( \frac1{n\pm it}\right) \left(
1+(n\pm it)^2 \right)^{(\n-1)/2}\, dt
\end{align*}
The integral with~$G_+$ is $\oh\left( \|\ps\|
n^{\re\n-1} \right) =o(1)$. Next we note that
\[ \frac 1t \sum_\pm \pm \ps\left( \frac1{n\pm
it}\right) \allowbreak \left( 1+(n\pm it)^2
\right)^{(\n-1)/2} = \oh\left( n^{\re\n-3} \|\ps'\| +
n^{\re\n-2} \|\ps\|\right) \]
for $0<t<1$. So the other integral is~$o(1)$ as well.

\ntst{Contribution of $\k_0(\n)$.}If $\n\in -2\NN$,
then $g_\k(\tau) = \oh(1)$ for $\re\tau \geq 2$. This
can be shown by the method of~\ref{kapnrepr} (deform
the line of integration into a narrow~$I_+$). So
$I_n(g_\k)= \oh(1/n)$.

For $\n\not\in -2\NN$ we take $g_\k(\tau) = \tau^{1+\n}
\allowbreak
\left(1+\tau^{-2}\right)^{(1+\n)/2}\allowbreak e^{-\p
i\n \sign(\im\tau)/2}$, see~\ref{kap0cont}.
\begin{align*} I_n(g_\k)
&= \frac i\p \int_0^\infty \sum_\pm
\ps\left(\frac1{n\pm it}\right) (n\pm it)^{\n-1}
\left( 1+(n\pm it)^{-2} \right)^{(1+\n)/2} e^{\mp \p i
\n/2}\, dt
\\
&= n^\n \frac1\p \int_0^\infty \sum_\pm \oh(\|\ps\|)
\left(t\mp i\right)^{\n-1} \left( 1+\oh(n^{-2}
\right)\, dt.
\end{align*}
The integral converges if $\re\n<0$. Under that
condition we find $I_n(g_\k)=\oh(n^{\re\n})=o(1)$.

\bw{Necessity}The integral $I_n(g_\k)$ with
$\ps(\tau)=\left(1+\tau^2\right)^{(1+\n)/2}$ equal
$n^\n$ times a non-zero function of~$\n$. So the bound
$\re\n<0$ is needed if $A_0(\al)\neq0$.

\ntst{Contribution of~$\m$.}\label{endaddprf}
$I_n(g_\m) = \lim_{T\rightarrow\infty} \frac 1\p
\int_{n-iT}^{n+iT} \ps\left(\frac1\tau\right)
\frac{-i}2 \tau \frac{d\tau}{1+\tau^2} = \frac12\ps(0)
$, for each~$n$.

\nwsect{Universal covering group}\label{seccov}
To see that the geodesic decomposition of automorphic
hyperfunctions is closely related to $1$-cocycles it
is better not to work on $\PSL_2(\RR)$, but on its
universal covering group.

\ntst{Universal covering group.} The universal covering
group~$\tilde G$ is a central extension of
$G=\PSL_2(\RR)$ with center $\tilde Z\cong \ZZ$. As an
analytic variety it is isomorphic to $\bhv\times\RR$.
This isomorphism is written as $(z,\th) \mapsto \tilde
p(z) \tilde k(\th)$ which covers the isomorphism
$\bhv\times\left( \RR\bmod\p\ZZ\right) \rightarrow G :
(z,\th) \mapsto p(z)k(\th)$.

There are injective continuous group homomorphisms
$\RR\rightarrow G:x\mapsto \tilde n(x)$,
$\RR^\ast_{>0}\rightarrow G:y\mapsto \tilde a(y)$, and
$\RR\rightarrow G:\th\mapsto \tilde k(\th)$, covering
respectively $x\mapsto \matr1x01$, $y\mapsto
\matc{\sqrt y}00{1/\sqrt y}$, and $\th \mapsto
k(\th)$. We have $\tilde p(z)=\tilde n(x)\tilde a(y)$.
The center of~$\tilde G$ is $\tilde Z\isdef \tilde
k(\p\ZZ)$. The group~$\tilde P \isdef \tilde p(\bhv)$
is isomorphic to~$P$; it is the connected component
of~$1$ in the parabolic subgroup $\tilde Z\tilde P$
of~$\tilde G$.

The projection $\tilde G\rightarrow G$ we write as
$g\mapsto \hat g$. We define a lifting $\SL_2(\RR)
\rightarrow\tilde G:g\mapsto \tilde g$ by
$\widetilde{\matr abcd} \tilde p(z) = \tilde
p\left(\frac{az+b}{cz+d} \right) k\left( -\arg(cz+d)
\right)$, with $-\p<\arg\leq\p$. Some properties:
$\widetilde{\matr abcd}^{-1} = \widetilde{\matr
d{-b}{-c}a}$ if $\arg\left(ci+d\right)\in (-\p,\p)$,
$\widetilde{p(z)}=\tilde p(z)$, and $\widetilde{k(\th)
} = \tilde k(\th)$ for $-\p<\th<\p$.

$\tilde \G$ is the full original in~$\tilde G$ of the
discrete subgroup~$\G$. It contains~$\tilde Z$.

\ntst{}\label{T}$\tilde T\isdef \tilde P \backslash
\tilde G$ is a covering of~$T=P\backslash G$. We use
the coordinate $\th$ corresponding to $\th\mapsto
\tilde P\tilde k(\th)$. The covering map $\tilde
T\rightarrow T$ corresponds to $\pr: \th \mapsto \tau=
\cot\th$. We denote the right action of $\tilde G$
on~$\tilde T$ in terms of the coordinate~$\th$ by
$g:\th\mapsto \th\cdot g$. We have $\th\cdot \tilde
k(\th_1)=\th+\th_1$. If $\arg\left(ci+d
\right)\in(-\p,\p)$, then $\th \mapsto
\th\cdot\widetilde{\matr abcd}$ is the strictly
increasing analytic function given by $\th\mapsto
\arg\left( (ia-b)\sin\th+(d-ic)\cos\th\right)$ on a
neighborhood of~$\th=0$. This function satisfies
$\left(\th+\p\right) \cdot\widetilde{\matr abcd} =
\th\cdot\widetilde{\matr abcd} +\p$.

The automorphism~$j$ of~$G$ is covered by the
automorphism $\tilde p(z)\tilde k(\th) \mapsto \tilde
p(-\bar z)\tilde k(-\th)$ of~$\tilde G$, also denoted
by~$j$. It induces an involution~$\spgl$ in the
functions on~$\tilde G$, which corresponds to the
reflection $\th \mapsto -\th$ in~$\tilde T$.

\ntst{Representations.}Let $\left( \tilde\p_\n,\tilde
M^\n\right)$ be the representation of~$\tilde G$ in
the functions on~$\tilde G$ that transform on the left
according to $\tilde p(z) \mapsto y^{(1+\n)/2}$. (Note
that the $\tilde Z$-behavior is not fixed. This is an
induced representation from~$\tilde P$ to~$\tilde G$.
Usually one would induce from $\tilde Z\tilde P$
to~$\tilde G$.) This representation can be realized in
the functions on~$\tilde T$. This leads to the action
$\tilde \p_\n(\tilde G)$ in the analytic functions
$\tilde M^\n_\om \isdef \Akr(\RR)$ on~$\tilde T$, and
in the space of hyperfunctions~$\tilde
M^\n_{-\om}\isdef \Bkr(\RR)$. We identify $\tilde T$
with~$\RR$ by means of the coordinate~$\th$. Note that
$\tilde\p_\n(\tilde k(\z))$ acts as the translation
over~$\z$.

Let $\tilde M^\n_b\isdef \Bkr_b(\RR)$ be the subspace
of $\tilde M^\n_{-\om}$ of hyperfunctions with bounded
support. It is invariant under $\tilde \p_\n$, and
there is a duality between $\left(\tilde \p_\n, \tilde
M^\n_\om \right)$ and $\left( \tilde \p_{-\n},
M^{-\n}_b \right)$, that can be described with
explicit integrals, see~\ref{dualR}.

Functions on~$G$ correspond to functions on~$\tilde G$
that are invariant under the center~$\tilde Z$. In
this way we obtain the following identifications:
\[ M^\n_\om = \left( \tilde M^\n_\om\right)^{\tilde Z},
\quad M^\n_{-\om} = \left( \tilde
M^\n_{-\om}\right)^{\tilde Z}, \quad A^\n_{-\om}(\G) =
\left( \tilde M^\n_{-\om}\right)^{\tilde \G}. \]
These identifications respect the reflection~$\spgl$.
The action of $\tilde \p_\n(g)$ in the $\p_\n(\tilde
Z)$-invariant spaces on the right corresponds to the
action of~$\p_\n(\hat g)$ in the spaces on the left.

\ntst{Projection.}There is also a natural map $\s:
\tilde M^\n_b \rightarrow M^\n_{-\om}$. If $\bt\in
\tilde M^\n_b$ has support inside an interval~$I$ of
length smaller than~$\p$ it corresponds, via
composition with $\p:\th\mapsto \cot\th$, to a
hyperfunction on~$T$ with support inside~$\pr(I)$.
This hyperfunction we define to be~$\s\bt$. We
extend~$\s$ to~$\tilde M^\n_b$ by additivity, using a
decomposition based on partings as discussed
in~\ref{decompR}. This linear map~$\s$ respects the
reflection~$\spgl$ and satisfies $\p_\n(\hat g) \s =
\s \tilde \p_\n(g)$. The kernel of~$\s$ consists of
the elements of the form $\left( \p_\n(\tilde
k(\p))-1\right) \bt_1$ with $\bt_1\in \tilde M^\n_b$.

\ntst{Cohomology.}Let $ X\subset\tilde T$ be invariant
under~$\tilde\G$. Let $A$ be a $\tilde\G$-module. We
consider the complex $C^\cdot_X(\tilde\G,A)$ defined
by
\begin{align*}
C^n_X(\tilde\G,A) &\isdef \vzm{c:X^{n+1}\rightarrow
A}{c(\x_0\cdot\g,\ldots,\x_n\cdot\g) = \g^{-1}
c(\x_0,\ldots,\x_n)}\\
dc(\x_0,\ldots,\x_{n+1})
&\isdef \sum_{l=0}^{n+1} (-1)^{l+1}
c(\x_0,\ldots,\widehat{\x_l},\ldots,\x_{n+1}).
\end{align*}
$Z^n_X(\tilde\G,A)$, $B^n_X(\tilde\G,A)$, and
$H^n_X(\tilde\G,A)$ are the corresponding groups of
cocycles, coboundaries and cohomology classes.

If $\G$ has only one cuspidal orbit, and $X$ is the
full original in~$\tilde T$ of the set of cusps, then
$H^1_X(\tilde\G, A)$ is the usual parabolic
cohomology, in which the inhomogeneous
$1$-cocycle~$\eta$ have the additional property
$\eta(\p) \in \left(\p-1\right) A$ if $\p$ fixes an
element of~$X$. (To a homogeneous $1$-cocycle~$c$
corresponds $\eta_c : \g\mapsto
c\left(0\cdot\g^{-1},0\right)$.)

\ntst{$\G$-decomposition and
$1$-cocycles.}\label{altoc}Suppose that $\al \in
A^\n_{-\om}(\G) = \left( \tilde
M^\n_{-\om}\right)^{\tilde\G}$ has a
$\G$-de\-com\-po\-si\-tion~$p$ on~$X$. Let $\tilde X$
be the full original of~$X$ in~$\tilde T$. We
construct a cocycle $c(\al,p) \in Z^1_{\tilde
X}(\tilde \G , \tilde M^\n_b)$ in the following way:
\begin{enumerate}
\item[a)] If $\x,\eta\in \tilde X$, $\x<\eta<\x+\p$,
then $c(\al,p;\eta,\x)$ is the hyperfunction with
support inside $[\x,\eta]_p$ that satisfies $\s
c(\al,p;\eta,\x) =\al[\pr(\eta),\pr(\x)]_p$.

The properties of $\al[\cdot,\cdot]_p$ imply
$c(\al,p;\eta,\th) + c(\al,p;\th,\x) =
c(\al,p;\eta,\x)$ if $\x<\th<\eta<\x+\p$, and
$c(\al,p;\eta\cdot\g,\x\cdot\g) = \p_\n (\g)^{-1}
\left( c(\al,p;\eta,\x) \right)$ for all
$\g\in\tilde\G$.

\item[b)] If $\eta \geq \x+\p$, then we take
intermediate points $\th_0=\x<\th_1<\cdot<
\th_k=\eta$, with $\th_j<\th_{j-1}+\p$, and define
$c(\al,p;\eta,\x) = \sum_{j=1}^k
c(\al,p;\th_{j},\th_{j-1})$. This does not depend on
the choice of the intermediate points.
\item[c)] $c(\al,p;\x,\x) \isdef 0$, $c(\al,p;\x,\eta)
\isdef -c(\al,p;\eta,\x)$ if $\eta<\x$.
\end{enumerate}
This turns out to define a $1$-cocycle. It has the
additional properties that
$\supp\left(c(\al,p;\x,\eta) \right) $ is contained in
the closed interval in~$\tilde T$ with end points $\x$
and~$\eta$, and that $c(\al,p;\eta,\x)$ and the
$\tilde\G$-invariant hyperfunction on~$\tilde T$
corresponding to~$\al$ have the same restriction to
the open interval $(\x,\eta)$.

If $p$ is the geodesic decomposition of~$\al$, we write
$c_\al$ instead of $c(\al,p)$.

If $j(\G)=\G$, then $\spgl(c_\al(\x,\eta)) =
c_{\spgl\al}(-\eta,-\x)$.

\begin{stel}{Proposition}\label{coh}Let $X$ be a
non-empty $\G$-invariant subset of~$T$. Denote its
full original in~$\tilde T$ by~$\tilde X$.
\begin{enumerate}
\item[i)] Let $\al \in A^\n_{-\om}(\G,X)$.
\begin{enumerate}
\item[a)] The class $[\al]$ of $c(\al,p)$ in
$H^1_{\tilde X}(\tilde \G, \tilde M^\n_b)$ does not
depend on the choice of the $\G$-decomposition~$p$.
\item[b)] There is a bijective map from the set of
$\G$-decompositions of~$\al$ on~$X$ to the set of
$h\in C^0_{\tilde X}(\tilde\G,\tilde M^\n_b)$ that
satisfy $\supp(h(\x))\subset\vz\x$ for all
$\x\in\tilde X$.
\end{enumerate}
\item[ii)] The map $\al \mapsto [\al]$ is an injection
$A^\n_{-\om}(\G,X) \rightarrow H^1_{\tilde X}(\tilde
\G, \tilde M^\n_b )$. The image consists of those
classes that have a representative $c\in C^1_{\tilde
X}(\tilde\G, \tilde M^\n_b)$ satisfying
$\supp(c(\x,\eta)) \subset [\x,\eta]$ for all
$\x,\eta\in\tilde X$, $\x<\eta$.
\end{enumerate}
\end{stel}
\bw{Remark}$A^\n_{-\om}(\G,X)$ is the subset of $\al
\in A^\n_{-\om}(\G)$ that have a $\G$-decomposition
on~$X$, see~\ref{Gamdecompdef}.
\bw{Proof} See the Lemmas~\ref{ctoal}--\ref{altocinj}.

\begin{stel}{Lemma}\label{ctoal}For each $c\in
C^1_{\tilde X}(\tilde\G, \tilde M^\n_b)$ that
satisfies $\supp(c(\x,\eta)) \subset [\x,\eta]$ for
all $\x,\eta\in\tilde X$, $\x<\eta$, there is a unique
$\al \in A^\n_{-\om}(\G, X)$ and a unique
$\G$-decomposition~$p$ of~$\al$ on~$X$ such that
$c=c(\al,p)$.
\end{stel}
\bw{Proof}Let $c$ be given. For $x,y\in T$, $x\neq y$,
we can find $\x,\eta\in\tilde X$ such that
$x=\pr(\x)$, $y=\pr(\eta)$, and $\eta<\x<\eta+\p$. We
define $A(x,y)\isdef \s c(\x,\eta)$. So $\supp(A(x,y))
\subset [x,y]$. The freedom in the choice of $\x$
and~$\eta$ is a translation over a multiple of~$\p$,
hence the choice of $\x$ and~$\eta$ does not influence
the definition of~$A(x,y)$. If there are $\al$ and~$p$
with $c(\al,p)=c$, then $\al[x,y]_p=A(x,y)$.

The properties $A(x,z)+A(z,y)=A(x,y)$ and $A(\g\cdot x,
\g\cdot y) =\p_\n(\g)A(x,y)$ are easily checked.

Now consider $x,y,z,u\in T$, $x\neq y$, $z\neq u$.
Choose $\x,\eta,\z,\ups\in\tilde T$ above these
elements such that $\eta<\x<\eta+\p$ and $\ups < \z <
\ups+\p$. Put $p_1\isdef c(\x,\eta) +
c\left(\eta,\x-\p\right)$, and $p_2 \isdef c(\z,\ups)
+ c\left(\ups,\z-\p\right)$. Then $A(x,y)+A(y,x) =\s
p_1$, and $A(z,u)+A(u,z) = \s p_2$. The cocycle
properties show that $\s\left(p_2-p_1\right) =
\s\left( c(\x,\z) + c\left( \z-\p, \x-\p \right)
\right) =0 $. This implies that $\al \isdef
A(x,y)+A(y,x)$ does not depend on the choice of $x\neq
y$. We have $\al \in A^\n_{-\om}(\G)$ and $\al[x,y]_p
\isdef A(x,y)$ defines a $\G$-decomposition of~$\al$
such that $c(\al,p)=c$.

\begin{stel}{Lemma}Let $\al \in A^\n_{-\om}(\G,X)$, and
let $p$ be a $\G$-decomposition of~$\al$ on~$X$. Then
each $\G$-decomposition~$q$ of~$\al$ on~$X$ has the
form~$q=p\langle h\rangle$, where
\[ \al[\pr(\x),\pr(\eta)]_{p\langle h\rangle} \isdef
\al[\pr(\x),\pr(\eta)]_p + \s h(\eta) -\s h(\x),\]
for $\eta<\x<\eta+\p$, $\x,\eta\in \tilde X$, where $h$
runs through the elements of $C^0_{\tilde X}(\tilde\G,
\tilde M^\n_b)$ that satisfy $\supp(h(\x)) \subset
\vz\x$ for $\x\in\tilde X$.

For such~$h$ we have $c(\al,p\langle h\rangle) =
c(\al,p) + dh$.
\end{stel}
\bw{Proof}The correspondence in Lemma~\ref{ctoal}
associates $\al$ and~$p\langle h\rangle$ to the
cocycle $c(\al,p)+dh$. This shows that $p\langle
h\rangle$ is a $\G$-decomposition.

Let $p$ and~$q$ be $\G$-decompositions of~$\al$ on~$X$.
One can check in~\ref{altoc} that the cocycle $c\isdef
c(\al,q) - c(\al,p)$ satisfies
$\supp(c(\x,\eta))\subset \vz{\x,\eta}$. So for
$\x\neq\eta$ we have $c(\x,\eta) = c_l(\x,\eta) +
c_r(\x,\eta)$ with $\supp(c_l(\x,\eta)) \subset\vz\x$
and $\supp(c_r(\x,\eta))\subset \vz\eta$. A
consideration of the cocycle relation for three
different points shows that $c_l(\x,\eta) = c_l(\x)$,
$c_r(\x,\eta) = c_r(\eta)$, and $c_l(\x) = - c_r(\x)$,
the $\tilde\G$-behavior is $c_r(\x\cdot \g) =
\tilde\p_\n(\g) c_r(\g)$. Take $h\isdef c_r$. Then
$c(\al,q) = c(\al,p)+dh$, with $h$ satisfying the
condition in the lemma.

\begin{stel}{Lemma}\label{altocinj}Let $\al \in
A^\n_{-\om}(\G,X)$ and let $p$ be a $\G$-decomposition
of~$\al$ on~$X$. If $c(\al,p) \in B^1_{\tilde
X}(\tilde \G, \tilde M^\n_b)$, then $\al=0$.
\end{stel}
\bw{Proof}Suppose $c(\al,p)=dh$ for some $h\in
C^0_{\tilde X}(\tilde\G,\tilde M^\n_b)$. Fix some
$\th\in\tilde X$. The support of~$h(\th)$ is contained
in some bounded closed interval~$J$. Take $N\in\NN$
large, such that $N\p + \min(J) > \max(J) + 4\p$. So
there is a closed interval~$I$ of length~$2\p$ between
$\th$ and $\th+N\p$ that does not intersect $J$ or
$J+N\p$.
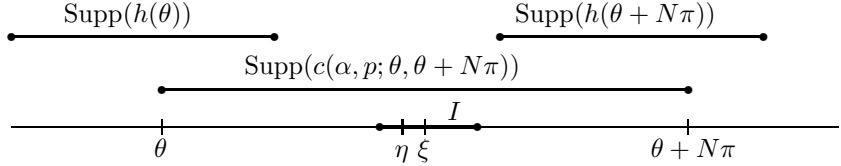
\begin{center}
\setlength\unitlength{1cm}
\begin{picture}(11,2)(0,-.4)
\put(0,0){\line(1,0){11}}
\put(2,-.1){\line(0,1){.2}}
\put(1.9,-.4){$\th$}
\put(9,-.1){\line(0,1){.2}}
\put(8.5,-.4){$\th+N\p$}
\put(5.2,-.1){\line(0,1){.2}}
\put(5.1,-.4){$\eta$}
\put(5.5,-.1){\line(0,1){.2}}
\put(5.4,-.4){$\x$}
\put(3.1,.7){$\supp(c(\al,p;\th,\th+N\p))$}
\put(.7,1.4){$\supp(h(\th))$}
\put(6.7,1.4){$\supp(h(\th+N\p))$}
\put(5.8,.1){$I$}
\thicklines
\put(0,1.2){\line(1,0){3.5}}
\put(0,1.2){\circle*{.1}}
\put(3.5,1.2){\circle*{.1}}
\put(6.5,1.2){\line(1,0){3.5}}
\put(6.5,1.2){\circle*{.1}}
\put(10,1.2){\circle*{.1}}
\put(2,.5){\line(1,0){7}}
\put(2,.5){\circle*{.1}}
\put(9,.5){\circle*{.1}}
\put(4.9,0){\line(1,0){1.3}}
\put(4.9,0){\circle*{.1}}
\put(6.2,0){\circle*{.1}}
\end{picture}
\end{center}
The restriction of $c\left(\al,p;\th,\th+N\p\right)$
to~$I$ vanishes. So for any $\x,\eta \in I$,
$\eta<\x<\eta+\p$, we have $\supp \left( c(\al,\p;
\x,\eta) \right) \subset \vz{\eta,\x}$, and hence
$\al[\pr(\x),\pr(\eta)]_p$ has support contained in
$\vz{\pr(\x),\pr(\eta)}$. So $\al$ has restriction
zero on each open interval in~$T$ bounded by two
different points of~$X$. As $X$ is infinite, we have
$\al=0$.

\nwsect{Image in a fixed weight}\label{secrho}
To return to Lewis's period function we make a
transition from hyperfunctions to functions. In
Proposition~\ref{cocycleintegral} we associate to
automorphic hyperfunctions cocycles with holomorphic
functions as values. Proposition~\ref{psiiotaalpha}
shows that this leads to solutions of~\vgl{MFE}.

\ntst{Functions of complex weight.}Let $q\in\CC$. We
define $C^\infty(\tilde G)_q$ to be the set of
functions in $C^\infty(\tilde G)$ that satisfy
$f(g\tilde k(\th))= f(g)e^{iq\th}$. This is the space
of functions of weight~$q$. These spaces are invariant
under left translation by elements of~$\tilde G$.

Functions in $f\in C^\infty(\tilde G)_q$ are fully
determined by the corresponding functions $z\mapsto
y^{-q/2} \allowbreak f(\tilde p(z))$ on~$\bhv$. The
left translation $L_{g_1} f(g) = f(g_1 g)$ in
$C^\infty(\tilde G)_q$ is a right $\tilde G$-action.
It corresponds to the right $\tilde G$-action $g_1 : F
\mapsto F|_q g_1$ in the functions on~$\bhv$ defined
by $ F|_q \tilde k(\p m) = e^{\p iqm} F$ for
$m\in\ZZ$, and $F|_q \widetilde { \matr abcd } (z)=
(cz+d)^{-q} F \left( \frac{az+b}{cz+d} \right)$ for
$\matr abcd \in \SL_2(\RR)$, $-\p<\arg(ci+d)<\p$.

\ntst{Map to weight~$q$.}The space~$\tilde M^{-\n}_\om$
contains $\ph:\th\mapsto e^{iq\th}$ for
each~$q\in\CC$. If we view $\tilde M^{-\n}_\om$ as a
space of functions on~$\tilde G$, then $\ph_q$ spans
the intersection $\tilde M^{-\n}_\om \cap
C^\infty(\tilde G)_q$.

If $q=2r\in2\ZZ$, then this function is invariant under
translation by~$\p$, and corresponds to the
function~$\ph_{2r}\in \Akr_T(T)$.

By $\r_q \bt : g \mapsto \left\langle
\tilde\p_{-\n}(g)\ph_q, \bt\right\rangle$ we define a
linear map $\tilde M^\n_b \rightarrow C^\infty(\tilde
G)_q$. It satisfies $\r_q \circ \tilde\p_\n(g) =
L_{g^{-1}} \circ \r_q$.

The map~$\r_q$ induces a homomorphism $H^1_{\tilde
X}(\tilde \G,\tilde M^\n_b) \rightarrow H^1_{\tilde
X}(\tilde\G, C^\infty(\tilde G) )$, with the left
action $\g\mapsto L_{\g^{-1}}$ of~~$\tilde\G$ on
$C^\infty(\tilde G)$. This may be uninteresting for
general weights~$q$, as the cohomology group with
values in $C^\infty(\tilde G)_q$ may vanish. (This is
the case for $\G=\Gmod$ and $q\not\in 2\ZZ$.)

\ntst{Image in weight~$1-\n$.}The differentiation
relations in~\ref{diffrel} extend to~$\ph_q$ with
arbitrary~$q$. This implies that $\EE^\pm \left(\r_q
\bt\right) = \left(1-\n\pm q\right)\r_{q\pm 2}\bt$.
For weight $q=1-\n$ we have $\EE^- \left(\r_{1-\n}
\bt\right)=0$. So the corresponding functions
on~$\bhv$ are holomorphic. We define the linear
map~$\R$ from $\tilde M^\n_b$ to the holomorphic
functions on~$\bhv$ by
\[ \R\bt (z) \isdef y^{(\n-1)/2} \left \langle \tilde
\p_{-\n}(\tilde p(z)) \ph_{1-\n} ,\bt\right\rangle. \]
$\R\bt$ is a holomorphic function on~$\bhv$, satisfying
\[ \left(\R\p_\n\left( \widetilde{\matr abcd} \right)
\right) \bt(z) = (a-cz)^{\n-1} (\R\bt)\left(
\frac{dz-b}{-cz+a}\right) = \left. \left(\R\bt\right)
\right|_{1-\n} \widetilde {\matr abcd}^{-1} (z) \]
for $\matr abcd\in \SL_2(\RR)$ with
$\arg\left(ci+d\right)\in
(-\p,\p)$.%
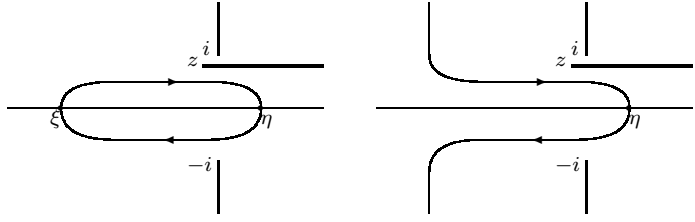
\begin{figure}
\begin{center}
\footnotesize
\setlength\unitlength{.7cm}
\begin{picture}(6,4)(-4,-2)
\put(-4,0){\line(1,0){6}}
\put(.8,0){\circle*{.1}}
\put(.8,-.3){$\eta$}
\put(-3.2,-.3){$\x$}
\put(-3,0){\circle*{.1}}
\qbezier(-3,0)(-3,.5)(-2,.5)
\qbezier(.8,0)(.8,.5)(-.2,.5)
\put(-2,.5){\line(1,0){1.8}}
\qbezier(-3,0)(-3,-.6)(-2,-.6)
\qbezier(.8,0)(.8,-.6)(-.2,-.6)
\put(-2,-.6){\line(1,0){1.8}}
\put(-1.05,.5){\vector(1,0){.3}}
\put(-.75,-.6){\vector(-1,0){.3}}
\put(-.3,1){$i$}
\put(-.6,-1.2){$-i$}
\put(-.6,.8){$z$}
\thicklines
\put(-.3,.8){\line(1,0){2.3}}
\put(0,1){\line(0,1){1}}
\put(0,-1){\line(0,-1){1}}
\end{picture}
\qquad
\begin{picture}(6,4)(-4,-2)
\put(-4,0){\line(1,0){6}}
\put(.8,0){\circle*{.1}}
\put(.8,-.3){$\eta$}
\put(-.3,1){$i$}
\put(-.6,-1.2){$-i$}
\put(-.6,.8){$z$}
\qbezier(.8,0)(.8,.5)(-.2,.5)
\qbezier(-2,.5)(-3,.5)(-3,1)
\put(-2,.5){\line(1,0){1.8}}
\put(-3,1){\line(0,1){1}}
\qbezier(.8,0)(.8,-.6)(-.2,-.6)
\qbezier(-2,-.6)(-3,-.6)(-3,-1.2)
\put(-.2,-.6){\line(-1,0){1.8}}
\put(-3,-1.2){\line(0,-1){.8}}
\put(-1.05,.5){\vector(1,0){.3}}
\put(-.75,-.6){\vector(-1,0){.3}}
\thicklines
\put(-.3,.8){\line(1,0){2.3}}
\put(0,1){\line(0,1){1}}
\put(0,-1){\line(0,-1){1}}
\end{picture}
\end{center}
\caption[o]{Contours used in the integral
representation of $\left(\R c_\al(\tilde\x,
\tilde\eta) \right)(z)$. In the picture on the right
$\x=\infty$.} \label{psiintcont}
\end{figure}

\begin{figure}
\begin{center}
\footnotesize
\setlength\unitlength{.7cm}
\begin{picture}(6,4)(-2,-2)
\put(-2,0){\line(1,0){6}}
\put(-.8,0){\circle*{.1}}
\put(-1.1,-.3){$\eta$}
\qbezier(-.8,0)(-.8,.5)(.2,.5)
\qbezier(2,.5)(3,.5)(3,1)
\put(2,.5){\line(-1,0){1.8}}
\put(3,1){\line(0,1){1}}
\qbezier(-.8,0)(-.8,-.6)(.2,-.6)
\qbezier(2,-.6)(3,-.6)(3,-1.2)
\put(.2,-.6){\line(1,0){1.8}}
\put(3,-1.2){\line(0,-1){.8}}
\put(.75,.5){\vector(1,0){.3}}
\put(1.05,-.6){\vector(-1,0){.3}}
\put(-.3,1){$i$}
\put(-.6,-1.2){$-i$}
\put(.5,.8){$z$}
\thicklines
\put(.4,.8){\line(-1,0){2.4}}
\put(0,1){\line(0,1){1}}
\put(0,-1){\line(0,-1){1}}
\end{picture}
\end{center}
\caption[o]{Contour used in the integral representation
of $\left(\R c_\al(\tilde\x, -\p) \right)(z)$. This
picture corresponds to the case $\eta=\infty$.}
\label{psiintcont1}
\end{figure}
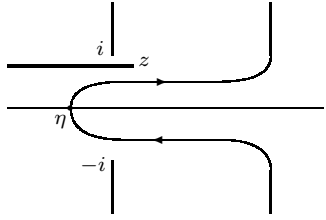

\begin{stel}{Proposition}\label{cocycleintegral}Suppose
that $\al\in A^\n_{-\om}(\G)$ has geodesic approach on
the elements of a $\G$-invariant set $X\subset T$. Let
$g$ be a representative of~$\al$.

Consider $\x,\eta\in X$, $\x\neq\eta$, such that
$\infty\not\in(\x,\eta)$. Choose
$\tilde\x>\tilde\eta\in \tilde T \cap [-\p,0]$ such
that $\pr\tilde\x=\x$ and $\pr\tilde\eta=\eta$. Then
\[ \left(\R c_\al(\tilde \x, \tilde \eta) \right)(z) =
\frac1\p \pvint_{Q(\x,\eta)} (z-\tau)^{\n-1}
\left(1+\tau^2\right)^{-(1+\n)/2} g(\tau) \, d\tau ,\]
with the contour $Q(\x,\eta)$ as indicated in
Figure~\ref{Qcont} on page~\pageref{Qcont}.

The function $z\mapsto \left(\R c_\al(\tilde\x,
\tilde\eta) \right)(z)$ has a holomorphic extension
from~$\bhv$ to $\CC\setminus (-\infty,\eta]$ if
$\eta\neq\infty$, and an extension to $\CC\setminus
[\x,\infty)$ if $\x\neq\infty$. If
$\infty<\x<\eta<\infty$, then the extension to~$\ohv$
across $(-\infty,\x)$ is $e^{2\p i\n}$ times the
extension across $(\eta,\infty)$.
\end{stel}
\bw{Remarks}In the Figures~\ref{psiintcont}
and~\ref{psiintcont1} we have drawn the contours and
the lines where $\arg\left(z-\tau\right)$ and
$\arg\left(1+\tau^2\right)$ jump. The standard choice
of the argument does not work if $\eta=\infty$. Then
 we take $\arg\left(z - \tau\right)\in(-2\p,0)$, see
 Figure~\ref{psiintcont1}. If we extend $\R
c_\al(\tilde\x,\tilde\eta)$ into~$\ohv$, we have to
adapt the choice of $\arg\left(z-\tau\right)$
continuously. In all cases the curve~$Q(\x,\eta)$
should stay inside the domain of the
representative~$g$ of~$\al$. The principal value
interpretation of the integral has been introduced
in~\ref{pv}.

If $\infty<\x<\eta<\infty$, then we can replace
$\left(1+\tau^2\right)^{-(1+\n)/2} g(\tau)$ by
$f_\al(\tau)$, see Proposition~\ref{alphtof}.

\bw{Proof}As $\x\neq\eta$, the length of
$[\tilde\eta,\tilde\x]$ is strictly less than~$\p$.
Hence the hyperfunction $c_\al(\tilde\x,\tilde\eta)$
on~$\tilde T$ corresponds to the hyperfunction
$\al[\x,\eta]$ on~$T$ with support inside $[\x,\eta]$.
We compute the quantity $\left\langle
\tilde\p_{-\n}(\tilde p(z))\ph_{1-\n},c_\al(\tilde \x,
\tilde \eta) \right\rangle $ as $\left\langle h_z,
\al[\x,\eta] \right\rangle$, where $h_z$ is the
holomorphic function on a neighborhood of~$[\x,\eta]$
such that $\th \mapsto h_z(\pr(\th))$ coincides with
$\th \mapsto \tilde\p_{-\n}(\tilde p(z)
)\ph_{1-\n}(\th)$ on a neighborhood of
$[\tilde\eta,\tilde\x]$.

We have $\tilde \p_{-\n}(\tilde p(z))\ph_{1-\n}(\th) =
y_1^{(1-\n)/2} e^{i(1-\n)\th_1}$, where $\tilde p(z_1)
\tilde k(\th_1) = \tilde k(\th) \tilde p(z)$. For
$-\p<\th <0$ we find $\tilde \p_{-\n}(\tilde
p(z))\ph_{1-\n}(\th) = y^{(1-\n)/2} \left(
\cos\th-z\sin\th\right)^{\n-1}$. For $\th\in (-\p,0)$
and $\tau=\cot\th \in T\setminus\vz\infty$ we have
$\arg\left( \cos\th -z\sin\th \right) = \arg
\frac{z-\tau}{\sqrt{1+\tau^2}} =
\arg\left(z-\tau\right)-\frac12\arg\left(1+\tau^2\right)$.
Hence $h_z(\tau) = y^{(1-\n)/2} \allowbreak \left(z
-\tau\right)^{\n-1} \allowbreak
\left(1+\tau^2\right)^{(1-\n)/2} $, on a neighborhood
of $[\x,\eta]$. This leads to the integral
representation in the proposition. The holomorphic
extensions into~$\ohv$ clearly exist.

Let $\infty<\x<\eta<\infty$. The continuation of $\R
c_\al(\tilde\x,\tilde\eta)$ across $(-\infty,\x)$ is
computed with the jump of $\arg\left(z-\tau\right)$ on
a curve that passes above $Q(\x,\eta)$, and the other
continuation with $\arg\left(z-\tau\right)$ jumping on
a curve below~$Q(\x,\eta)$. So in the former integral
the argument of $z-\tau$ is $2\p$ more than in the
other.%

\ntst{Reflection.}\label{jc}If we choose originals
$\hat\x$ and~$\hat\eta$ in $[0,\p]$, then we get a
similar integral for $\R c_\al(\hat\x,\hat\eta)$, with
$\left(z-\tau\right)^{\n-1}$ replaced by
$\left(\tau-z\right)^{\n-1}$.

Let $j(\G)=\G$, and $\eta<\infty$. If we use
$\tau\mapsto -g(-\tau)$ as the representative
of~$\spgl\al$, we obtain $\left(\R \spgl\left(
c_\al(\tilde\x,\tilde\eta) \right) \right)(z) =
\left(\R c_\al(\tilde\x,\tilde\eta)\right)(-z)$, where
we use the continuation of $\R
c_\al(\tilde\x,\tilde\eta)$ across $(\eta,\infty)$.

\ntst{Holomorphic cusp form.}If $\al\in
A^{2k-1}_{-\om}(\G)$ corresponds to a holomorphic cusp
form of weight~$2k$, then $\R
c_\al(\tilde\x,\tilde\eta) (z) $ is a polynomial
in~$z$ that turns out to be $(-1)^k 4^{1-k}$ times the
period polynomial $R_H(\x,\eta)$ introduced
in~\ref{clcocycle}.

\ntst{Eisenstein series.}\label{Bzf}Take $\x=\infty$
and $\eta=0$. {}From~\ref{Eisdecomp} we find for $\re
s>1$ and $z\in\CC\setminus(-\infty,0)$:
\begin{align*}
&\left( \R c_{\e_s}\left(0,-\frac\p2 \right) \right)(z)
\\
&\qquad= \frac12\p^{-s} (s-1) \sideset{}'\sum_{p,q}
f_{p,q} \left(p^2+q^2\right)^{-s} \left. \left(
(z-\tau)^{-2s} (1+\tau^2)^s \right)
\right|_{\tau=-p/q} \\
&\qquad= \p^{-s}(s-1) \left( \frac12\z(2s)(1+z^{-2s}) +
\sum_{p,q\geq1} (qz+p)^{-2s} \right),
\end{align*}
with $f_{p,q}$ as in~\ref{Eisdecomp}, and $\arg(z),\,
\arg\left(qz+p\right) \in (-\p,\p)$.

\ntst{Interpretation.}Let $\al$ be as in
Proposition~\ref{cocycleintegral}. This proposition
shows that for $\tilde \x, \tilde\eta$ projecting
to~$\G\cdot\infty$ with distance strictly smaller
than~$\p$ the image $\R c_\al(\tilde\x,\tilde\eta)$
has an extension into~$\ohv$ across part of~$T$. In
this way we can view $ \R c_\al$ as a cocycle with
values in the $\tilde G$-module of holomorphic
functions with domain of the form $\prc\setminus I$,
where $I\subset T$ depends on the function. The action
of $\widetilde{\matr abcd}$ with
$\arg\left(ci+d\right) \in (-\p,\p)$ is given by the
same formula as above, with the standard choice of the
argument of $cz+d$. This forces $F|_q \tilde k(\p) =
e^{\pm\p iq} F$ on~$\hv^\pm$.
\medskip

\ntst{Modular group.}\label{modgrcoh}The covering
$\tilde\Gmod$ of~$\Gmod$ is generated by $w\isdef
\tilde k\left(\frac\p2\right)$ and $n\isdef \tilde
n(1)$, subject to the relations $w^2n=nw^2$ and
$nwnwn=w^3$.

Let $ X = \prq$, and $\tilde X=\pr^{-1}X$. Any $c\in
C^1_{\tilde X}(\tilde\Gmod,A)$ for a
$\tilde\Gmod$-module~$A$ is determined by
$p=c\left(0\cdot w^{-1},0\right) = c\left(
-\frac\p2,0\right)$, with the relations
$\left(nwn+n\right) p = \left(1+w+w^2\right)p$ and
$\left(n-1\right) \left(w+1\right)p=0$. (To check
this, it is convenient to work with the corresponding
inhomogeneous cocycle $\g\mapsto c\left(0\cdot
\g^{-1},0\right)$.)

The cocycle is a coboundary if and only if
$p=\left(1-w\right)f$ for some $f\in A$ satisfying
$nf=f$. If we can solve $f $ from $\left(1+w\right)p =
\left(1-w^2\right) f$, then $c\in B^1_{\tilde
X}(\tilde\Gmod,A)$. This equation is always solvable
if $A$ is the module of holomorphic functions on
$\bhv\cup\ohv$ with the action of weight~$q
\in\CC\setminus2\ZZ $ indicated above. This is used in
Lemma~\ref{FkrtoPsikr}.

If $\al$ satisfies the assumptions of
Proposition~\ref{cocycleintegral}, and $c=c_\al$, then
$p=\R c_\al\left(-\frac\p2,0\right)$. The relation
$c_\al\left(0,-\frac\p2\right) = c_\al\left(
0,-\frac\p4\right) + c_\al \left( -\frac\p4, -\frac\p
2\right)$ corresponds to $p = n^{-1}p + wnw^{-1}p$.
This is the functional equation for the $\ps$-function
in~\vgl{MFEnu}, but with weight $1+\n$ replaced by
$1-\n$. In that case
$F\isdef\left(1-w^2\right)f=\left(1+w\right)p $ is
given by
\begin{align*}
F(z) &= \frac{-1}\p \left( \pvint_{Q(\infty,0)} +
\pvint_{Q(0,\infty)} \right) (z-\tau)^{\n-1}
\left(1+\tau^2\right)^{-(1+\n)/2} g(\tau)\, d\tau \\
&= \frac{-1}\p \pvint_{I_+-I_-} (z-\tau)^{\n-1}
\left(1+\tau^2\right)^{-(1+\n)/2} g(\tau)\, d\tau,
\end{align*}
with $I_\pm$ as in Figure~\ref{Icont} on
page~\pageref{Icont}, and $\arg\left(z-\tau\right) \in
\left(-\frac\p2,\frac32\p\right)$ if $z\in\bhv$ and in
$\left(-\frac32\p,\frac\p2 \right)$ if $z\in\ohv$. By
$I_+-I_-$ we mean the union of both contours, with
opposite orientation of~$I_-$. The function~$F$ on
$\bhv\cup\ohv$ is holomorphic and has period~$1$. It
can be defined for other~$\G$ as well.

\begin{stel}{Proposition}\label{psiiotaalpha}Suppose
that $\al\in A^\n_{-\om}(\G)$ has polynomial growth
at~$\infty$ and satisfies one of the conditions a)
and~b) in Theorem~\ref{thmgdc}, implying geodesic
decomposition on~$\G\cdot\infty$. Define
\[ F_\al(z) \isdef \frac{-1}\p \pvint_{I_+-I_-}
(z-\tau)^{\n-1} \left(1+\tau^2\right)^{-(1+\n)/2}
g(\tau)\, d\tau, \]
where $g$ represents~$\al$, and $I_\pm$ and
$\arg\left(z-\tau\right)$ are as indicated above.
Then, for $z\in\hv^\pm$,
\begin{align*}
F_\al(z) &= \pm \frac{2i e^{\pm\p i\n/2}}{ (2\p)^\n
\Gf(1-\n) } \sum_{m=1}^\infty m^{-\n} A_{\pm m}(\al)
e^{\pm 2\p imz}
\\
&\quad \hbox{} \pm ie^{\pm \p i\n/2} \sin\frac{\p\n}2
B_0(\al) \pm \text{(if $\n=0$) } iA_0(\al).
\end{align*}

If $\n\not\in 1+2\ZZ$, and $\n\not\in 2\NN$, then
$F_\al(z)= 2i\sqrt\p e^{\pm\p i\n/2} \allowbreak
\Gf\left( \frac{1+\n}2\right)^{-1} \allowbreak
\Gf\left( 1-\frac\n2\right)^{-1} \allowbreak
f_{\iota(\n)\al}$ on~$\hv^\pm$, with
$f_{\iota(\n)\al}$ the function associated to
$\iota(\n)\al \in A^{-\n}_{-\om}(\G)$ in
Proposition~\ref{alphtof}.

If $\n\not\in 1+2\ZZ$, $\n\not\in 2\NN$, and
$\G=\Gmod$, then $\R c_\al\left(-\frac\p2,0 \right) =
\frac i{\sqrt\p} \Gf\left( \frac{1-\n}2\right)
\allowbreak \Gf\left(1 -\frac\n2\right)^{-1}
\allowbreak \ps_{\iota(\n)\al}$, where
$\ps_{\iota(\n)\al} \in \Psimod(-\n) $ is the image of
$\iota(\n)\al$ under the map in
Theorem~\ref{authftopsi}.
\end{stel}
\bw{Remark}See \vgl{falphdef} and
Lemma~\ref{cuspdecomp} for the coefficients $A_n(\al)$
and~$B_0(\al)$, and \vgl{iotadef} on
page~\pageref{iotadef} for the definition of the
isomorphism~$\iota(\n)$.
\bw{Proof}For any hyperfunction~$\bt$ with
representative~$h$ we define
\[ J_\n(z,h) = \frac{-1}\p \pvint_{I_+-I_-}
(z-\tau)^{\n-1}
\left(1+\tau^2\right)^{-(1+\n)/2}h(\tau)\,d\tau.\]
This does not depend on the choice of the
representative~$h$, so we also write $J_\n(z,\bt)$. We
have $F_\al(z) = J_\n(z,\al)$.

Take $z\in\hv^\z$, with $\z\in \vz{1,-1}$. By taking
the contours~$I_\pm$ wide, we see that
$J_\n\left(z-x,h\right) = J_\n\left( z,
\p_\n\matr1{-x}01 h\right)$. This implies that
$F_\al(z) = J_\n(z,\al)$ has a Fourier expansion
$F_\al(z) = \sum_{m\in\ZZ} F(m,\z) e^{2\p imz}$ for
$z\in\hv^\z$, and $F(m,z)$ is given by $ e^{-2\p imz}
J_\n(z,\Fkr_m\al)$, see~\ref{Ftal}.

In~\ref{genappr}--\ref{kapbyint} we determine
$J_\n(z,h)$ for representatives~$h$ of the various
Fourier $N$-equi\-variant hyperfunctions that can
occur in~$\Fkr_m\al$. These computations give the
Fourier expansion of~$F_\al$ indicated in the
proposition.

Let $\n\in\CC\setminus\left(1+2\ZZ\right)$. The
relation with $\iota(\n)\al$ follows from~\ref{isoAB}.
In the modular case we have already seen that
$F_\al(\tau) = p(\tau) + \tau^{\n-1} p(-1/\tau)$, with
$p=\R c_\al\left(-\frac\p2,0\right)$. The relations in
Lemma~\ref{FkrtoPsikr} allow us to express
$\ps_{\iota(\n)\al}$ in terms of~$p$.

\ntst{General remarks.}\label{genappr}Let $N>|\re x|$
and $ 0<\e < |\im z|$. We have $J_\n= - J_\n^0 -
\sum_{\dt=1,\,-1} \dt J_\n^\dt$, with
\begin{align*}
J_\n^0(z,h) &= \sum_\pm \frac{\pm1}\p \int_{\x=-N}^N
\left( z-\x\mp i\e\right)^{\n-1} \left( 1+(\x\pm
i\e)^2\right)^{-(1+\n)/2} h(\x\pm i\e)\,d\x,\\
J_\n^\dt(z,h) &= \frac i\p \int_{t=\e}^\infty \sum_\pm
(z-\dt N \mp it)^{\n-1} \left(1+(\dt N\pm it)^2
\right)^{-(1+\n)/2} h(\dt N\pm it )\, dt\\
&= \frac i\p \int_{t=\e}^\infty \sum_\pm f_\dt(\dt N\pm
it) \frac{ h(\dt N\pm it ) } { 1+(\dt N\pm it)^2
}\,dt,
\end{align*}
with $f_\dt(\tau) = -\dt e^{\p i\z(1+\dt)\n/2}
\left(1-\frac z\tau\right)^{\n-1} \left( 1+\tau^{-2}
\right)^{(1-\n)/2}$.

\ntst{Case $n=0$, hyperfunction $\m$.}For the
representative $h(\tau) = \frac{-i}2 \tau$
representing~$\m$ we can take $\e=0$. Then $J^0_\n$
vanishes, and $J_\n^1(z,h) = \left\langle f_1, \m[
N,\infty] \right \rangle = \frac12 f_1(\infty) = -
\frac12 e^{\p i\z\n} $, and $J_\n^{-1}(z,h) =
-\left\langle f_{-1}, \m[\infty,-N] \right\rangle =
-\frac12 f_{-1}(\infty) = -\frac12$. Hence $J_\n(z,\m)
= \z i e^{\p i\z\n/2} \sin\frac{\p \n}2$.

\ntst{Case $m=0$, hyperfunction $\k_0(\n)$.}We consider
$\n\in\CC\setminus(-2\NN)$, with $\re\n<1$. Take $p$
and~$q$ as in~\ref{kap0cont}. We find
\begin{align*}
J_\n^\dt(z,q) &= \frac{\dt e^{\p i\z(1+\dt)\n/2}
\Gf(1+\n/2) }{2\p\sqrt\p} \int_\e^\infty \sum_\pm \mp
i \left( t \mp i\dt N \pm iz \right)^{\n-1} \, dt
\\
&= \frac{\dt i e^{\p i\z(1+\dt)\n/2}
\Gf(\n/2)}{4\p\sqrt \p} \sum_\pm \pm \left( \e \mp
i\dt N \pm iz \right)^\n,\\
J_\n^0(z,p) &= \frac{-1}{2\sqrt\p \Gf(1-\n/2)} \left(
e^{\p i\z\n} (-z+N+i\e)^\n - (z+N-i\e)^\n\right).
\end{align*}
The limits for $\e\downarrow0$ exist, and after taking
$\e=0$ it does no longer hurt that we use different
representatives. As $N\rightarrow\infty$, all terms
tend to zero, and we obtain $J_\n(z,\k_0(\n)) = 0 $ if
$\n\not\in -2\NN$. (One can also check for finite~$N$
that the terms cancel each other.)

\ntst{Case $m=0$, $\n=0$, hyperfunction $\ld(0)$.}We
find \[ J_0(z,\ld(0) ) = \lim_{\n\rightarrow0}
\frac2{\n\sqrt\p} \left( J_\n(z,\k_0(\n) ) -
\frac1{\sqrt\p} J_\n(z,\m) \right) = -\z i.\]

\ntst{$\k_m(\n)$ given by an
integral.}\label{kapbyint}Let $\re \n<0$ if $m=0$. We
write a representative of~$\k_m(\n)$ as
in~\ref{kapnrepr}:
\[ g(\tau) = \frac1{2\p i} \int_{\tilde I}
\frac{1+\tau\tau_0}{\tau_0-\tau} e^{2\p im\tau_0}
\left(1+\tau_0^2\right)^{(\n-1)/2}\, d\tau_0,\]
where $\tilde I$ is a path of the form $I_+$ if
$m\geq0$, and of the form~$I_-$ otherwise. In both
cases we take $\tilde I$ between the contours $I_+$
and~$I_-$ that we use to compute $J_\n(z,\k_m(\n))$.
If we insert this representation of~$g$ into the
definition of $J_\n(z,\k_m(\n))$, everything converges
absolutely, provided we take the principal value
interpretation of the outer integral. Hence we can
interchange the order of integration, and find:
\begin{align*}
J_\n(z,\k_m(\n)) &= \frac1{2\p i} \int_{\tilde I}
\frac{-1}\p \pvint_{I_+-I_-} (z-\tau)^{\n-1}
\left(1+\tau^2\right)^{-(1+\n)/2}
\frac{1+\tau\tau_0}{\tau_0-\tau} \,d\tau\\
&\qquad\qquad\hbox{} \cdot e^{2\p im\tau_0}
\left(1+\tau_0^2\right)^{(\n-1)/2}\, d\tau_0.
\end{align*}
Next we enlarge the contours~$I_\pm$ into wider
contours~$\hat I_\pm$ of the same type, but such that
$\tilde I$ is contained in one of them. This changes
the inner integral into
\[ \pvint_{\hat I_+-\hat I_-} (z-\tau)^{\n-1}
\left(1+\tau^2\right)^{-(1+\n)/2}
\frac{1+\tau\tau_0}{\tau_0-\tau} \,d\tau + 2\p i
(z-\tau_0)^{\n-1} \left(1+\tau_0^2\right)^{(1-\n)/2}.
\]
The integrand in the integral over~$\hat I_\pm$ is
$\oh\left( \tau^{-2}\right)$ as
$|\tau|\rightarrow\infty$. So we do not need a
principal value interpretation. The integral can be
deformed into integrals over vertical lines at $\re
\tau=-M$ and $\re\tau=M$. The limit as $M$ tends to
infinity yields~$0$. We are left with
\[ J_\n(z,\k_m(\n)) = \frac{-1}\p \int_{\tilde I}
(z-\tau_0)^{\n-1} e^{2\p im\tau_0}\,d\tau_0.\]

If $m=0$, this vanishes by an explicit computation. If
$\sign \im z \neq \sign m$, then we move off the
path~$\tilde I$ upwards or downwards, and obtain
$J_\n(z,\k_m(\n))=0$ if $\z m<0$. In the case $\z m>0$
we have a holomorphic function of~$\n$ that we compute
for $\re\n$ large, and find $J_\n(z,\k_m(\n)) = 2i\z
e^{\p i\z\n/2} \allowbreak \Gf\left(1-\n\right)^{-1}
\allowbreak (2\p|m|)^{-\n} \allowbreak e^{2\p imz}
$.\medskip

\ntst{Discussion.}For $\al\in A^\n_{-\om}(\Gmod)$
corresponding to a cuspidal Maass form, and for
$\al=\e_s$ with $\re s<1$, $s\not\in\ZZ$, the period
function~$\ps_\al$ has turned out to arise in two
different ways:
\begin{itemize}
\item Theorem~\ref{authftopsi} and its proof show that
an analysis of representatives of~$\al$ leads
to~$\ps_\al$. (This works for all $\al \in
A^\n_{-\om}(\G)$.)
\item The geodesic decomposition of $\iota(\n)\al$
leads to a cocycle on~$\tilde\G$ with values in the
hyperfunctions on~$\tilde T$ with bounded support.
Testing against the lowest weight vector~$\ph_{1+\n}$
in~$\tilde M^\n_\om$ leads to a cocycle with values in
the holomorphic functions on $\bhv\cup\ohv$. This
cocycle is determined by one value, that turns out to
be a multiple of~$\ps_\al$.
\end{itemize}

\nwsect{Transfer operator}\label{sectrop}
We conclude this paper with some remarks on the
transfer operator. We restrict ourselves to results
that are a direct consequence of the previous
sections.

\ntst{Transfer operator.}\label{trop}Let $\al \in
A^\n_{-\om}(\Gmod)$ satisfy the assumptions of
Theorem~\ref{add}, and $B_0(\al)=0$. This means that
$\al$ corresponds to a cuspidal Maass form. Then
$\al[\infty,0] = \sum_{n=0}^\infty \al\left[
-n-1,-n\right]$ (weak convergence). This implies that
\[\R c_\al\left(-\frac\p2,0\right)(z) =
\sum_{n=0}^\infty \R c_\al\left( -\frac\p2+ \arctan n,
-\frac\p2 + \arctan (n+1) \right) (z)\]
for each $z\in\CC\setminus (-\infty,0)$. Put $p_\al =
\R c_\al\left(-\frac\p2,0\right)$. As
\[\spgl \tilde\p_\n\widetilde{\matc n {-n-1} 1{-1}
}\left( c_{\spgl\al}\left(-\frac\p2, 0\right) \right)
= c_\al\left( -\frac\p2+\arctan n, -\frac\p2 +
\arctan(n+1)\right),\]
we have \[\R c_\al\left( -\frac\p2+ \arctan n,
-\frac\p2 + \arctan (n+1) \right) (z) =
\left(z+n\right)^{\n-1} p_{\spgl\al}
\left(1+\frac1{n+z}\right)\]
(see~\ref{T} for the action on~$\tilde T$, and~\ref{jc}
for the reflection). So for $\spgl\al = \pm \al$, the
function $p_\al$, and its multiple
$\ps_{\iota(\n)\al}$, are eigenfunctions of the
transfer operator
 \[L_{1-\n} : f \mapsto \sum_{n=0}^\infty
 \left(z+n\right)^{\n-1} f\left(1+\frac1{n+z} \right)\]
 \ of Mayer, see~\cite{May}. (We have shifted the
 functions.) Note that the convergence implies $f(1)=0$
 if $\re\n\geq0$. This implies that $\ps_\al(1)=0$ for
$\al$ associated to cuspidal Maass forms.

\ntst{Transfer operator for hyperfunctions.}We define
\[\Lkr_\n \isdef \sum_{n=0}^\infty \spgl \p_\n\matc
n{-n-1}1{-1}\] as an operator in~$M^\n_{-\om}$. Its
domain consists of those $\bt \in M^\n_{-\om}$ with
support contained in $[\infty,0]$ for which the sum
converges weakly.

We have seen that, for $\al$ as above, $\al[0,1]$ is in
the domain of~$\Lkr_\n$, and is an eigenvector with
eigenvalue~$\pm1$ if $\spgl\al = \pm \al$.

\ntst{Eisenstein series.}Theorem~\ref{add} can be
applied to~$\e_s$ for $\re s>\frac12$. For these
values of~$s$ the hyperfunction $\e_s[\infty,0]$ is in
the domain of~$\Lkr_{1-2s}$, and
$\Lkr_{1-2s}\e_s[\infty,0] = \e_s[\infty,0]
-\frac12\p^{-s}\left(s-1\right) \z(2s)\m$. (We have
used $\spgl\e_s=\e_s$.)

This suggests that $\e_s[\infty,0]$ would be an
eigenvector of~$\Lkr_{1-2s}$ if $\z(2s)=0$. But we
have no continuation of our results to these values
of~$s$. After continuation of the expression for $\R
c_{\e_s}\left(0,-\frac\p2 \right)$ in~\ref{Bzf}, we
obtain an eigenfunction of~$L_{2s}$.

\newcommand\bibit[5]{
\bibitem
{#2}#3: {\em #4;\/ } #5}

\bigskip

\noindent R.W.Bruggeman, Mathematisch Instituut,
Universiteit Utrecht, Postbus~80.010, 3508 TA Utrecht,
The Netherlands

\noindent E-mail: {\tt bruggeman@@math.ruu.nl}

\end{document}